# AN RKHS FORMULATION OF THE INVERSE REGRESSION DIMENSION-REDUCTION PROBLEM

By Tailen Hsing[1] and Haobo Ren

*University of Michigan and Sanofi-Aventis Pharmaceutical*

Suppose that $Y$ is a scalar and $X$ is a second-order stochastic process, where $Y$ and $X$ are conditionally independent given the random variables $\xi_1, \ldots, \xi_p$ which belong to the closed span $L_X^2$ of $X$. This paper investigates a unified framework for the inverse regression dimension-reduction problem. It is found that the identification of $L_X^2$ with the reproducing kernel Hilbert space of $X$ provides a platform for a seamless extension from the finite- to infinite-dimensional settings. It also facilitates convenient computational algorithms that can be applied to a variety of models.

**1. Introduction.** Identifying the space spanned by the inverse regression function leads to a highly effective dimension-reduction approach for nonparametric regression function estimation. See Duan and Li (1991), Li (1991), Chen and Li (1998), Cook (1998) and Cook and Li (2002), among others. In this paper, we consider the approach in the context where the predictor is a stochastic process. Our goal is to introduce a unified formulation that can be applied to a wide variety of models, and, at the same time, retains the spirit of multivariate analysis so that statistical inference can be carried out in a natural and efficient manner.

Let $(\Omega, \mathcal{F}, P)$ be a probability space, and let $L^2(\Omega, \mathcal{F}, P)$ be the Hilbert space containing all random variables on $(\Omega, \mathcal{F}, P)$ that have finite variances, and with inner product defined by $\langle U, V \rangle_{L^2(\Omega, \mathcal{F}, P)} = E(UV)$. Let $Y$ be a random element defined on $(\Omega, \mathcal{F}, P)$. The nature of $Y$ critically influences the construction of the computational algorithms, but is of no relevance in the theoretical formulation of the inverse regression problem. Let $\{X_t, t \in T\}$ be a real-valued, zero-mean, second-order stochastic process defined on $(\Omega, \mathcal{F}, P)$, where the index set $T$ is assumed to be a separable metric space.

Received February 2007.
[1]Supported by NSF Grants DMS-06-24239 and DMS-07-07021.
*AMS 2000 subject classifications.* Primary 62H99; secondary 62M99.
*Key words and phrases.* Functional data analysis, sliced inverse regression.

Here, $T$ may be quite flexible which can be a single homogeneous set or a union of sets with different topological nature; for example, $T = \bigcup_{q=1}^{Q} T_q$, where $T_1 = [a,b]$, $T_2 = \{t_1, \ldots, t_J\}$, and so on, in which case one can think of the restrictions of $X_t$ to the $T_q$ as covariates of different functional nature. Note that we do not assume that the paths of $X_t$ lie in a known Hilbert space, which is a common assumption in functional data analysis literature [cf. Ramsay and Silverman (2005), Dauxois, Ferré and Yao (2001) and Ferré and Yao (2003, 2005)]. Indeed, in the infinite-dimensional case, such an assumption may be restrictive and the identification of the Hilbert space may pose an extra problem in practice. Eubank and Hsing (2007) contains a discussion on the theoretical limitations of this assumption.

As usual, the Hilbert space $L_X^2$ of $\{X_t, t \in T\}$ is defined as the subspace of $L^2(\Omega, \mathcal{F}, P)$ that contains all finite linear combinations of the form $\sum_{i=1}^{k} c_i X_{t_i}$, $t_i \in T, c_i \in \mathbb{R}, k = 1, 2, \ldots$ and their limits in $L^2(\Omega, \mathcal{F}, P)$. See Ash and Gardner (1975) for details of these notions.

Define the following conditions (IR1) and (IR2) in which $\xi_1, \ldots, \xi_p$ are fixed elements in $L_X^2$:

(IR1) $Y$ and $X$ are conditionally independent given $\xi_1, \ldots, \xi_p$.
(IR2) For any $\xi \in L_X^2$, $E(\xi | \xi_1, \ldots, \xi_p) \in \text{span}\{\xi_1, \ldots, \xi_p\}$ a.s.

A particularly relevant model for which (IR1) holds is the multiple-index model

$$(1) \qquad Y = \ell(\xi_1, \ldots, \xi_p, \varepsilon),$$

where $\varepsilon$ is a random error independent of the process $\{X_t\}$, and we call each $\xi_i$ an index and $\ell$ the link function. The number of indices, the indices themselves and the link function are all assumed unknown in practice.

Condition (IR2) holds if the joint distribution of any finite collection of elements from $L_X^2$ is elliptically contoured, which would be the case if, for instance, $\{X_t\}$ is a Gaussian process. However, this could be much more general [see Hall and Li (1993)]. It is clear that the indices $\xi_i$'s in (1) are nonidentifiable if $\ell$ is not specified. However, the $L^2$ subspace

$$L_{X,e}^2 := \text{span}\{\xi_1, \ldots, \xi_p\}$$

is identifiable. Following Li (1991), call $L_{X,e}^2$ the effective dimension-reduction space (EDRS) for (1). We are interested in estimating the EDRS, and in some situations, the link $\ell$.

It might be awkward to conceptualize the estimation of $L_{X,e}^2$ directly since it is a space of random variables. In some cases, this problem can be overcome naturally. For instance, if the sample paths of $X_t$ are contained in a Hilbert space $\mathcal{H}$ and $\xi_j = \langle \beta_j, X \rangle_\mathcal{H}$, where $\beta_j$ is the representer of the functional, then the problem of estimating $L_{X,e}^2$ can be solved by estimating the space spanned by the $\beta_j$. Indeed this is the approach adopted for the multivariate case in Li (1991) and for the functional data case in Ferré and Yao (2003, 2005). See also Dauxois, Ferré and Yao (2001). However, as mentioned earlier, we do not assume that the sample paths of $X_t$ are contained in a Hilbert space. Thus, we are interested in a natural and flexible representation of the $\xi_j$. Our solution is the reproducing kernel Hilbert space (RKHS) of $X_t$. It is known that the RKHS of $X_t$ is a mirror image of $L_X^2$ in terms of Hilbert space structure (cf. Section 2), and so the estimation of the EDRS in $L_X^2$ can, in principle, be accomplished through the estimation of the corresponding space in the RKHS. The primary goal of this paper is to show how this idea can be implemented, and the advantages of the approach. It is interesting to note that the possibility of such an RKHS formulation was mentioned very briefly in Remark 2.4 of Li (1992).

The structure of this paper is as follows. We review the basic properties of RKHS in Section 2. We do so out of the concern that the notion of RKHS is not a part of the standard statistics curriculum today and our readers may not be familiar with the relevant facts required in this paper. Section 3 contains a key theoretical result on the inverse regression function $\mathrm{E}(X_t|Y)$ that facilitates dimension-reduction. Estimation issues are addressed in Section 4, where an asymptotic theory will also be developed; the inference will be conducted based on the data $(x_i, y_i)$, $i = 1, \ldots, n$, with each $x_i$ observed at a finite set of points. In Section 5, we provide a number of numerical examples, including simulation studies and a data analysis. Finally, the proofs are collected in Section 6.

We should mention that the present paper focuses on the basic RKHS formulation of the inverse regression dimension-reduction problem but ignores many important theoretical and methodological aspects that go along with the formulation, such as tests for determining $p$ in (1), choice of the number of slices in the sliced inverse regression procedure, estimating smooth representers $\beta_j$ when $\xi_j = \langle \beta_j, X \rangle_{\mathcal{H}}$, and so on. They will hopefully be pursued in future works by those that find this approach meaningful.

**2. Reproducing kernel Hilbert spaces.** Since the seminal work of Parzen (1959, 1961a, 1961b, 1963), statistical innovations using RKHS have been steadily developed. See Wahba (1990), Gu (2002) and Berlinet and Thomas-Agnan (2004). A quick survey reveals that the notion of RKHS is now embraced strongly by the machine learning community due to its importance in regularization problems. In this section we present the general definitions and common properties of RKHS required in this paper. The details of most of the results can be found in Aronszajn (1950). Other relevant references will be provided in due course. In order to be self-contained, short proofs are provided whenever suitable in Section 6.

A symmetric, real-valued bivariate function $K$ defined on $T$ is said to be nonnegative definite, denoted by $K \geq 0$, if for all $n \in \mathbb{N}$, $a_1, \ldots, a_n \in \mathbb{R}$, and $t_1, \ldots, t_n \in T$, we have $\sum_{i,j=1}^n a_i a_j K(t_i, t_j) \geq 0$. For convenience, symmetric nonnegative definite bivariate functions will be referred to as covariance kernels below. Also, for any bivariate function $K$, write

$$K_t = K(\cdot, t).$$

DEFINITION 1. A Hilbert space $\mathcal{H}$ is said to be a RKHS if the elements of $\mathcal{H}$ are functions defined on some set $T$, and there is a bivariate function $K$ on $T \times T$, having the following two properties:

(a) For all $t \in T$, $K_t \in \mathcal{H}$.
(b) For all $t \in T$ and $f \in \mathcal{H}$, $f(t) = \langle f, K_t \rangle_\mathcal{H}$.

In this case, $K$ is said to be a reproducing kernel of $\mathcal{H}$.

The following fundamental result is known as the Moore–Aronszajn theorem.

PROPOSITION 1. (a) *If $K$ is a reproducing kernel of $\mathcal{H}$, then $K$ is a covariance kernel and is unique. Conversely, if $K$ is a covariance kernel on $T \times T$, a unique RHKS of functions on $T$ with $K$ as the reproducing kernel can be constructed.*
(b) *If $K$ is the reproducing kernel of the RKHS $\mathcal{H}$, then $\mathrm{span}\{K_t, t \in T\}$ is dense in $\mathcal{H}$.*

Property (b) of Definition 1, called the reproducing property, is the essence of the notion of RKHS and will be applied extensively throughout this paper. The notation $\mathcal{H}_K$ will be used to denote the RKHS having the reproducing kernel $K$.

An important reason why RKHS plays an important role in statistics is that the Hilbert space of a second-order stochastic process can be represented by the RKHS whose reproducing kernel equals the covariance function of the process. To see that, consider a second-order, zero-mean process $\{X_t, t \in T\}$ with covariance function $R$. As usual, $\mathcal{H}_R$ denotes the RKHS with reproducing kernel $R$. Consider the linear map $\Psi_X$ from $L_X^2$ to $\mathcal{H}_R$ satisfying

$$\Psi_X(X_t) = R_t, \qquad t \in T.$$

PROPOSITION 2. $\Psi_X$ *is an isometric isomorphism, namely, it is one-to-one and satisfies* $\langle \eta, \xi \rangle_{L_X^2} = \langle \Psi_X(\eta), \Psi_X(\xi) \rangle_{\mathcal{H}_R}, \eta, \xi \in L_X^2$.

The mapping $\Psi_X$ was introduced by Loève (1948) and is sometimes referred to as Loève's isometry. For more information on the duality between a stochastic process and its RKHS [see Wahba (1990)].

The following result given in Theorem 1.1 of Fortet (1973) provides an insightful way to compute the RKHS norm.

PROPOSITION 3. *A function $f$ on $T$ is in $\mathcal{H}_K$ iff*

$$\sup_{t_1,\ldots,t_n} \sup_{a_i} \frac{|\sum_{i=1}^n a_i f(t_i)|^2}{\sum_{i=1}^n \sum_{j=1}^n a_i a_j K(t_i, t_j)} < \infty, \tag{2}$$

*where the suprema are taken over all $t_1, \ldots, t_n \in T$ and all real $a_1, \ldots, a_n$ for all $n$, such that the denominator in (2) is nonzero. If $f \in \mathcal{H}_K$, then the left-hand side of (2) is the RKHS norm.*

PROPOSITION 4. *Suppose that $K_1$ and $K_2$ are two covariance kernels on $T \times T$ with $K_2 - K_1 \geq 0$. Then:*
  (a) $\mathcal{H}_{K_2} \supset \mathcal{H}_{K_1}$ *where* $\|f\|_{\mathcal{H}_{K_2}} \leq \|f\|_{\mathcal{H}_{K_1}}$ *for $f \in \mathcal{H}_{K_1}$, and,*
  (b) *the linear operator $L: \mathcal{H}_{K_2} \mapsto \mathcal{H}_{K_1}$ for which*

$$LK_2(\cdot, t) = K_1(\cdot, t), \qquad t \in T$$

*is a bounded, nonnegative definite, and self-adjoint operator on $\mathcal{H}_{K_2}$.*

DEFINITION 2. Under the assumption of Proposition 4, we say that $K_2$ dominates $K_1$ if $K_2 - K_1 \geq 0$, denoted by $K_2 \geq K_1$, and call $L$ the dominance operator of $\mathcal{H}_{K_2}$ over $\mathcal{H}_{K_1}$. If $L$ is nuclear, or trace-class, namely $L$ satisfies $\operatorname{tr}(L) < \infty$, we say that $K_2$ nuclear-dominates $K_1$, denoted by $K_2 \gg K_1$, and $L$ is called a nuclear dominance operator.

The trivial case when $K_2 = K_1$ can be provided as an illustration, where $L$ is the identity mapping. Whether $K_2 \gg K_1$ in this case, of course, depends on the dimensionality of $T$.

Let $T$ be an index set and $T_1 \subset T$. For any $f$ defined on $T$, let $f|_{T_1}$ stand for the restriction of $f$ to the subset of $T_1$.

PROPOSITION 5. *Let $T$ be a separable metric space of which $S_0 = \{s_1, s_2, \ldots\}$ is a dense subset. Let $K$ be a covariance kernel on $T \times T$ and $K_n = K|_{S_n \times S_n}$, where $S_n = \{s_1, \ldots, s_n\}$. For any function $f$ defined on $T$, write $f_n = f|_{S_n}$. The following hold:*
  (a) *For any function $f$ defined on $T$, if for some $n > 1$, $f_n \in \mathcal{H}_{K_n}$, then $f_m \in \mathcal{H}_{K_m}$ for any $m \leq n$ and*

$$\|f_m\|_{\mathcal{H}_{K_m}} \leq \|f_n\|_{\mathcal{H}_{K_n}}.$$

  (b) *Let $f_n \in \mathcal{H}_{K_n}$ for any $n$, and $\lim_{n \to \infty} \|f_n\|_{\mathcal{H}_{K_n}} < \infty$. If either $T$ is countable or both $K$ and $f$ are continuous functions defined on $T \times T$ and $T$, respectively, then $f \in \mathcal{H}_K$ and $\|f\|_{\mathcal{H}_K} = \lim_{n \to \infty} \|f_n\|_{\mathcal{H}_{K_n}}$.*

**3. The covariance operator of inverse regression.** Below we continue to use the notation developed in Sections 1 and 2, and assume that (IR1) and (IR2) hold. As in Section 1, $L^2_{X,e} = \text{span}\{\xi_1, \ldots, \xi_p\}$ denotes the EDRS of (1) in $L^2_X$. Define the counterpart of the EDRS in $\mathcal{H}_R$:

$$\mathcal{H}_{X,e} = \Psi_X(L^2_{X,e}) = \text{span}\{\Psi_X(\xi_1), \ldots, \Psi_X(\xi_p)\},$$

which we call the reproducing kernel EDRS. We wish to conduct inference on $\mathcal{H}_{X,e}$ and $L^2_{X,e}$.

Denote by $Z_t$ the inverse regression process $\text{E}(X_t|Y), t \in T$. Clearly, $Z_t$ is also a second-order stochastic process with mean 0. Denote the covariance function of $Z_t$ by $K$. For $m = 1, 2, \ldots, t_1, \ldots, t_m \in T$, let $\mathbf{X} = (X_{t_1}, \ldots, X_{t_m})^T$, $\mathbf{Z} = (Z_{t_1}, \ldots, Z_{t_m})^T$ and $\mathbf{a} = (a_1, \ldots, a_m)^T \in \mathbb{R}^m$. We have

(3) $$\text{var}(\mathbf{a}^T \mathbf{X}) = \text{var}(\mathbf{a}^T \mathbf{Z}) + \text{E}(\text{var}(\mathbf{a}^T \mathbf{X}|Y)).$$

This implies that

(4) $$R \geq K$$

and it follows from Proposition 4 that

$$\mathcal{H}_K \subseteq \mathcal{H}_R.$$

Motivated by Theorem 3.1 of Li (1991), we make the following claim:

(5) The sample paths of $Z_t$ are in $\mathcal{H}_{X,e}$ a.s.

We will establish the validity of (5) in Theorem 6 below. However, let us first assume that (5) holds and consider some implications. Since (5) implies that the sample paths of $Z_t$ are in $\mathcal{H}_R$ a.s., we can define the covariance operator

(6) $$L = \text{E}\left(Z \bigotimes_{\mathcal{H}_R} Z\right),$$

where the tensor product $g \otimes_{\mathcal{H}_R} h$ denotes the linear operator that maps $f$ to $\langle g, f \rangle_{\mathcal{H}_R} \cdot h$ for $f, g, h \in \mathcal{H}_R$. By the reproducing property,

$$(LR_t)(s) = \text{E}(\langle Z, R_t \rangle_{\mathcal{H}_R} Z_s) = \text{E}(Z_t Z_s) = K_t(s), \qquad s, t \in T,$$

which implies that $\text{Im}(L) = \mathcal{H}_K$ and $L$ is the dominance operator of $\mathcal{H}_R$ over $\mathcal{H}_K$ [cf. Definition 2 and (b) of Proposition 4]. On the other hand, it follows readily from (5) that $K_t \in \mathcal{H}_{X,e}$ for all $t \in T$ and hence

$$\text{Im}(L) = \mathcal{H}_K \subseteq \mathcal{H}_{X,e}.$$

Thus, $\dim(\mathcal{H}_K) \leq \dim(\mathcal{H}_{X,e}) = p$. In particular, if $\mathcal{H}_K = \mathcal{H}_{X,e}$, estimating the eigenfunctions of $L$ provides an approach for estimating $\mathcal{H}_{X,e}$. Clearly, establishing (5) is crucial.

REMARKS. (a) Note that $\mathcal{H}_K$ is not always equal to $\mathcal{H}_{X,e}$. See Cook (1998) for a thorough discussion on this and related issues. Extending the ideas in this paper for dealing with those situations will be a topic of future research.

(b) Let $X$ be multivariate, that is, finite-dimensional, and we denote it by $\mathbf{X}$ for clarity. As before, assume that $\mathbf{X}$ has mean 0, covariance matrix $\mathbf{R}$, and let $\mathbf{Z} = \mathrm{E}(\mathbf{X}|Y)$. Then $\mathcal{H}_\mathbf{R}$ contains elements spanned by the column vectors of $\mathbf{R}$, where

$$\langle \mathbf{f}, \mathbf{g} \rangle_{\mathcal{H}_\mathbf{R}} = \mathrm{E}(\mathbf{f}^T \mathbf{R}^- \mathbf{X}\mathbf{X}^T \mathbf{R}^- \mathbf{g}) = \mathbf{f}^T \mathbf{R}^- \mathbf{g} \tag{7}$$

and

$$\Psi_\mathbf{X}^{-1} \mathbf{f} = (\mathbf{R}^- \mathbf{f})^T \mathbf{X}, \mathbf{f} \in \mathcal{H}_\mathbf{R}, \tag{8}$$

$\mathbf{R}^-$ being the Moore–Penrose generalized inverse of $\mathbf{R}$. Li (1991) showed that $\mathbf{Z} \in L^2_{\mathbf{X},e}$ with probability 1, and it follows that

$$L\mathbf{g} = \mathrm{E}\left(\mathbf{Z} \bigotimes_{\mathcal{H}_\mathbf{R}} \mathbf{Z}\right) \mathbf{g} = \mathrm{E}(\langle \mathbf{Z}, \mathbf{g} \rangle_{\mathcal{H}_\mathbf{R}} \mathbf{Z}) = \mathrm{E}(\mathbf{Z}\mathbf{Z}^T) \mathbf{R}^- \mathbf{g} =: \mathbf{K} \mathbf{R}^- \mathbf{g}. \tag{9}$$

By (7) and (9) the eigenvectors of $L$ in $\mathcal{H}_\mathbf{R}$ are $\mathbf{R}^{1/2}\mathbf{g}_i$ with the $\mathbf{g}_i$ denoting the eigenvectors of $\mathbf{R}^{-1/2}\mathbf{K}\mathbf{R}^{-1/2}$ in the Euclidean space. Provided that $\mathcal{H}_\mathbf{K} = \mathcal{H}_{\mathbf{X},e}$, it follows from (8) that $L^2_{\mathbf{X},e}$ is estimated by the span of $\Psi_\mathbf{X}^{-1}(\mathbf{R}^{1/2}\mathbf{g}_i) = (\mathbf{R}^{-1/2}\mathbf{g}_i)^T \mathbf{X}$. This completely agrees with the result for the multivariate setting described in Li (1991).

At first glance, (5) might seem a straightforward extension of similar results for the multivariate case in Li (1991), or the functional case in Ferré and Yao (2003, 2005). However, a closer inspection reveals that this is not the case, and, to establish it, a deeper understanding of the relationship between a RKHS and the sample paths of a stochastic process is called for. To give an idea of where the difficulties lie, recall that Parzen (1963) showed that almost all the sample paths of $X$ lie outside of $\mathcal{H}_R$ if $T$ is an infinite separable metric space and $R$ is continuous on $T \times T$; for instance, if $X$ is standard Brownian motion on $[0,1]$, then the paths are nowhere differentiable with probability 1 but $\mathcal{H}_R$ contains functions with square-integrable derivatives. In those situations, $Z_t$ is an average of paths that are a.s. not in $\mathcal{H}_R$, let alone $\mathcal{H}_{X,e}$. Driscoll (1973) gave sufficient conditions under which the sample paths of a Gaussian process fall into a RKSH; Lukić and Beder (2001) provided a much more general treatment of this class of problems, going beyond Gaussianity. The following development is partly inspired by their results.

In addition to the conditions (IR1) and (IR2) in Section 1, define the following condition:

(IR3) Either $T$ is countable, or both $R$ is continuous on $T \times T$ and the sample paths of $\mathrm{E}(X_t|Y)$ are continuous on $T$ with probability 1.

The following can be proved.

THEOREM 6.  *Assume the conditions* (IR1)–(IR3). *Then (5) holds.*

The proof of Theorem 6 will be given in Section 5. An approach for estimating $L^2_{X,e}$ by estimating $L$ and $\Psi_X$ as explained above will be developed in Section 4.

To introduce sliced inverse regression (SIR) in $\mathcal{H}_R$, consider the stochastic process

$$Z_t^{\mathcal{G}} = \mathrm{E}(Z_t|\mathcal{G})$$

for a given $\sigma$-field $\mathcal{G}$. An example of $\mathcal{G}$ is the $\sigma$-field generated by the sets $\{\omega \in \Omega : Y \in I_s\}$, $s = 1, \ldots, S$, where the $I_s$, called slices in Li (1991), are disjoint sets forming a partition of the range of $Y$. Denote by $K^{\mathcal{G}}$ the covariance function of $Z_t^{\mathcal{G}}$. The same variance decomposition argument in (3) shows that $K \geq K^{\mathcal{G}}$, and Proposition 4 implies that $\mathcal{H}_{K^{\mathcal{G}}} \subseteq \mathcal{H}_K$. If the conditions of Theorem 6 hold, then we have

$$\mathcal{H}_{K^{\mathcal{G}}} \subseteq \mathcal{H}_K \subseteq \mathcal{H}_{X,e}. \tag{10}$$

Denote the dominance operator of $\mathcal{H}_R$ over $\mathcal{H}_{K^{\mathcal{G}}}$ by $L^{\mathcal{G}}$, which is the covariance operator

$$L^{\mathcal{G}} = \mathrm{E}\left(Z^{\mathcal{G}} \bigotimes_{\mathcal{H}_R} Z^{\mathcal{G}}\right). \tag{11}$$

As before, estimating the eigenfunctions of $L^{\mathcal{G}}$ also estimates $\mathcal{H}_{X,e}$ if $\mathcal{H}_{K^{\mathcal{G}}} = \mathcal{H}_{X,e}$.

**4. Estimation and asymptotic theory.**  Assume without further reference in this section that the conditions (IR1)–(IR3) hold so that the conclusion of Theorem 6 holds. The primary goals of this section are to describe a procedure of estimation based on SIR, and to develop an asymptotic theory for the procedure.

In view of the description in Section 3, the estimations of the covariance function $R$ and inverse regression covariance function $K$ are clearly crucial elements in this problem. In some cases, these could be done more efficiently if the precise nature of the sample paths of $X$ is known. For example, in the infinite-dimensional case if the sample paths of $X$ are $m$-times continuously differentiable for some $m$, the incorporation of the information in nonparametric estimation procedures may lead to a faster rate of convergence in estimating $R$ and $K$ [see Rice and Silverman (1991), Silverman (1996), James, Hastie and Sugar (2000), Ramsay and Silverman (2005) and Wu and Pourahmadi (2003)]. However, we will not make such assumptions here, as our aim is to consider a general procedure whose principles and properties will, for a large part, transcend the detailed nature of the path properties of the second-order stochastic process $X$. Indeed, the development below simultaneously addresses both the finite- and infinite-dimensional cases.

We continue to use the notation defined in Section 3. In addition, for a real symmetric, nonnegative-definite matrix $\mathbf{A}$, let $\lambda_j(\mathbf{A})$ be the $j$th largest eigenvalue of $\mathbf{A}$; if the eigendecomposition of $\mathbf{A}$ is

$$\mathbf{A} = \sum \lambda_j(\mathbf{A}) \mathbf{u}_j \mathbf{u}_j^T,$$

define the generalized power

$$\mathbf{A}^\alpha = \sum_{\lambda_j(\mathbf{A}) > 0} \lambda_j^\alpha(\mathbf{A}) \mathbf{u}_j \mathbf{u}_j^T, \qquad \alpha \in \mathbb{R}.$$

Note that $\mathbf{A}^{-1}$ is the Moore–Penrose inverse of $\mathbf{A}$ and we denote it by $\mathbf{A}^-$.

We will focus on estimating $L^2_{X,e}$ by SIR, that is, assume that $\mathcal{H}_{K^\mathcal{G}} = \mathcal{H}_{X,e}$, where $\mathcal{G}$ is the $\sigma$-field generated by the sets $\{\omega \in \Omega : Y \in I_s\}, s = 1, \ldots, S$, the sets $I_1, \ldots, I_S$ forming a partition of the range of $Y$ with

$$p_s := \mathrm{P}(Y \in I_s) > 0 \qquad \text{for each } s.$$

Note that $Z^\mathcal{G} = \mathrm{E}(\mathrm{E}(X|Y)|\mathcal{G}) = \mathrm{E}(X|\mathcal{G})$, so

$$L^\mathcal{G} = \mathrm{E}\bigg(\mathrm{E}(X|\mathcal{G}) \bigotimes_{\mathcal{H}_R} \mathrm{E}(X|\mathcal{G})\bigg) = \sum_{s=1}^{S} p_s h_s \bigotimes_{\mathcal{H}_R} h_s,$$

where

(12) $$h_s = \mathrm{E}(X|Y \in I_s).$$

To fix ideas, let the eigenvalues of $L^\mathcal{G}$ be distinct; for $1 \leq j \leq p$, let $f_j$ denote the eigenfunction corresponding to the $j$th largest eigenvalues of $L^\mathcal{G}$, and, without loss of generality, let $\xi_j = \Psi_X^{-1}(f_j)$.

Let $(Y_i, X_{i,t}), 1 \leq i \leq n$, be $n$ i.i.d. realizations of $(Y, X_t)$. However, we only observe $Y_i, X_{i,t_j}, 1 \leq i \leq n, 1 \leq j \leq J_n$, for some finite $J_n$. Let

$$\mathbf{X}_i = (X_{i,t_1}, \ldots, X_{i,t_{J_n}})^T \quad \text{and} \quad \mathbf{h}_s = \mathrm{E}(\mathbf{X}_1 | Y_1 \in I_s)$$

and, for each $J$,

(13) $$\mathbf{R}_J = \{\mathrm{E}(X_{1,t_i} X_{1,t_j})\}_{i,j=1}^J.$$

Estimate $p_s$, $\mathbf{h}_s$ and $\mathbf{R}_{J_n}$, respectively, by the empirical estimators

$$\widehat{p}_s = \frac{1}{n} \sum_{i=1}^{n} I(Y_i \in I_s), \qquad \widehat{\mathbf{h}}_s = \frac{\sum_{i=1}^{n} \mathbf{X}_i I(Y_i \in I_s)}{\sum_{i=1}^{n} I(Y_i \in I_s)} \quad \text{and} \quad \widehat{\mathbf{R}}_{n, J_n} = \frac{1}{n} \sum_{i=1}^{n} \mathbf{X}_i \mathbf{X}_i^T.$$

If $X_i$ is not centered, then we need to center where appropriate in $\mathbf{h}_s$ and $\widehat{\mathbf{R}}_{n,J_n}$. As mentioned in the beginning of this section, in practice, if $X$ has smooth paths, then incorporating the information in the estimation of $\mathbf{h}_s$ and $\mathbf{R}_J$ may lead to estimators that are more efficient than the naive ones defined here.

For $k \leq J$, let $\mathbf{P}_{J,k}$ and $\widehat{\mathbf{P}}_{n,J_n,k}$ be the projection matrices onto the eigenspaces of the first $k$ eigenvalues of $\mathbf{R}_J$ and $\widehat{\mathbf{R}}_{n,J_n}$, respectively; let

$$(14) \quad \mathbf{R}_{J,k} = \mathbf{P}_{J,k} \mathbf{R}_J \mathbf{P}_{J,k} \quad \text{and} \quad \widehat{\mathbf{R}}_{n,J_n,k} = \widehat{\mathbf{P}}_{n,J_n,k} \widehat{\mathbf{R}}_{n,J_n} \widehat{\mathbf{P}}_{n,J_n,k}.$$

Our proposed estimator of $\xi_j$ is

$$(15) \quad \widehat{\xi}_{n,k,j} = (X_{t_1}, \ldots, X_{t_{J_n}}) \widehat{\mathbf{R}}_{n,J_n,k}^{-1/2} \mathbf{v}_j =: (X_{t_1}, \ldots, X_{t_{J_n}}) \widehat{\beta}_{n,k,j},$$

where $\mathbf{v}_j$ is the eigenvector corresponding to the $j$th largest eigenvalue of

$$(16) \quad \widehat{\mathbf{M}}_{n,k} := \widehat{\mathbf{R}}_{n,J_n,k}^{-1/2} \left( \sum_{s=1}^{S} \widehat{p}_s \widehat{\mathbf{h}}_s \widehat{\mathbf{h}}_s^T \right) \widehat{\mathbf{R}}_{n,J_n,k}^{-1/2}$$

in $\mathbb{R}^{J_n}$. Note that $X$ in (15) is a generic process whose sole purpose is to facilitate the definition of the estimator $\widehat{\xi}_{n,k,j}$. Below we will investigate the convergence of $\widehat{\xi}_{n,k,j}$ to $\xi_j$ in $L_X^2$.

REMARKS. (a) The procedure described above is not entirely new. If $X_i$ is finite-dimensional with $J_n = J$, then taking $k = J$ reduces the procedure above to that in Li (1991). In the infinite-dimensional case, the estimation of the eigenspaces of $\mathbf{R}_{J_n}$ corresponding to small eigenvalues is typically unstable, in which case $k$ acts as a smoothing parameter that controls the trade-off between bias and variance. Chiaromonte and Martinelli (2002) considered a similar approach in the context of analyzing gene-expression data.

(b) Ferré and Yao (2003) assume that the paths of $X_i$ are in a known Hilbert space $\mathcal{H}$. Their procedure is a "continuous" version of ours, since they assume that functional data $X_i$ are observed in their entirety. Of course, functional data are never observed in their entirety, so some kind of discrete approximation will have to be incorporated to implement their procedure. As such, there is little difference between their procedure and ours in that setting.

(c) In the infinite-dimensional case, if the observational points are different for different $X_i$, then smoothing of observed data $\mathbf{X}_i$ becomes necessary. In that case, the quantities $\widehat{\mathbf{h}}_s$ and $\widehat{\mathbf{R}}_{n,J_n,k}^{-1/2}$ will be computed based on the smoothed data. The details of this will be worked out in future work.

We proceed to explain the motivations of $\widehat{\xi}_{n,k,j}$ and develop an asymptotic theory. Let $J_n$ be nondecreasing and tending to some $J_\infty$ as $n \to \infty$, where $J_\infty$ is assumed to be $\infty$ for the infinite-dimensional case. A related issue for the infinite-dimensional case is that we stated earlier that we observe $X_{i,t_j}, 1 \leq i \leq n, 1 \leq j \leq J_n$, at stage $n$, but we did not specify the manner in which the set of observation points $t_1, \ldots, t_{J_n}$ change with $n$. There are two options in that regard. The first one is to consider the fully general case where each $t_j$ actually also depends on $n$ so that $t_j = t_{n,j}$. Another option is to consider a nested sequence of sets

$$T_J := \{t_1, \ldots, t_J\}, \qquad J \geq 1,$$

where the $t_j$ do not depend on $J$, so that more observation points will simply be added to each $X_i$ as $n$ increases. As far as the proofs go, the two cases require similar arguments. However, since the nested-sequence assumption entails slightly simpler details and much cleaner notation, we will take that approach.

First define two technical conditions, both of which amount to requiring that the leading eigenvalues of $\mathbf{R}_J$ dominate the rest. The first condition is

(17) $$\lim_{k \to k_\infty} \limsup_{J \to J_\infty} \operatorname{tr}((\mathbf{R}_J^- - \mathbf{R}_{J,k}^-)\mathbf{K}_J) = 0 \qquad \text{for some } k_\infty \leq J_\infty,$$

where $\mathbf{K}_J = \{K(t_i, t_j)\}_{i,j=1}^J$, and $\mathbf{R}_J$ and $\mathbf{R}_{J,k}$ are as defined in (13) and (14), respectively. It is shown by Lemma 12 below that, under very general conditions, $\lim_{J \to J_\infty} \operatorname{tr}(\mathbf{R}_J^- \mathbf{K}_J) < \infty$, which can be shown to be the trace of the dominance operator from $\mathcal{H}_R$ to $\mathcal{H}_K$. In that light, the condition (17) is quite mild.

To motivate the second technical condition, observe that if $0 < \inf_t \operatorname{E}(X_t^2) \leq \sup_t \operatorname{E}(X_t^2) < \infty$, then

(18) $$\operatorname{tr}(\mathbf{R}_J) = \sum_{j \geq 1} \lambda_j(\mathbf{R}_J) = \sum_{j=1}^J \operatorname{E}(X_{t_j}^2) = O(J).$$

If the random variables $X_{t_i}, 1 \leq i \leq J$, are uncorrelated, $\mathbf{R}_J$ is a diagonal matrix and all of the eigenvalues are bounded away from 0 and $\infty$. In the infinite-dimensional case, we wish to avoid this type of situation and focus on those where the strength of dependence among the $X_{t_i}$ increases as $J$ increases, so that the leading eigenvectors of $\mathbf{R}_J$ dominate. In that case, gaps of size $O(J)$ can be expected to exist between leading eigenvalues. The second technical condition is, for $m$ equal to a fixed positive integer,

(19) $$\liminf_{J \to J_\infty} \frac{\rho_m(\mathbf{R}_J)}{J} > 0,$$

where

(20) $$\rho_m(\mathbf{R}_J) = \min\{|\lambda_j(\mathbf{R}_J) - \lambda_m(\mathbf{R}_J)| : \lambda_j(\mathbf{R}_J) \neq \lambda_m(\mathbf{R}_J)\}.$$

Indeed, the conditions (17) and (19) are extremely general, as reflected by the following result.

PROPOSITION 7. *Let $T = [a, b]$ be any compact interval, $t, \ldots, t_J$ be equally spaced in $T$. If the covariance function $R$ is continuous on $T \times T$, then (17) holds with $k_\infty = J_\infty = \infty$. If, additionally, the multiplicity of $\lambda_m(Q)$ is 1, where $Q$ is the integral operator $Q : f \to \int_a^b R(\cdot, y) f(y) \, dy, f \in L^2[a, b]$, then (19) holds as well.*

The first step in establishing the estimator (15) is to compare $L^\mathcal{G}$ with the following operator:

$$\widetilde{L}_{n,k}^\mathcal{G} = \sum_{s=1}^{S} \widehat{p}_s \widetilde{h}_s \bigotimes_{\mathcal{H}_R} \widetilde{h}_s,$$

where

(21) $$\widetilde{h}_s = \widetilde{h}_{s,n,k} = (R(\cdot, t_1), \ldots, R(\cdot, t_{J_n})) \mathbf{R}_{J_n,k}^{-} \widehat{\mathbf{h}}_s.$$

Let $\|\cdot\|_\infty$ denote the sup or uniform norm of an operator.

LEMMA 8. *Assume that either $T = \bigcup_{J=1}^\infty T_J$, or $\bigcup_{J=1}^\infty T_J$ is dense in $T$ and $R$ is continuous on $T \times T$. Also assume that (17) holds. Then we have*

$$\|\widetilde{L}_{n,k}^\mathcal{G} - L^\mathcal{G}\|_\infty \xrightarrow{p} 0 \qquad \text{as } n \to \infty \text{ and then } k \to k_\infty.$$

Under the conclusion of Lemma 8, the eigenvalues and eigenfunctions of $\widetilde{L}_{n,k}^\mathcal{G}$ converge in probability to those of $L^\mathcal{G}$, where convergence of the eigenfunctions is in terms of the norm of $\mathcal{H}_R$. The convergence of the eigenvalues follows from Corollary 4 on page 1090 of Dunford and Schwarz (1988). The convergence of the eigenfunctions follows from the convergence of projection operators of eigenspaces, which can be established as in Gohberg and Kreĭn (1969), page 15 [see also Dauxois, Pousse and Romain (1982), pages 141–142].

Next, we express the eigenproblem of $\widetilde{L}_{n,k}^\mathcal{G}$ as an eigenproblem in $\mathbb{R}^J$.

LEMMA 9. *Let $(\lambda_j, \mathbf{u}_j)$ be the eigenvalues and eigenvectors of*

(22) $$\widetilde{\mathbf{M}}_{n,k} := \mathbf{R}_{J_n,k}^{-1/2} \left( \sum_{s=1}^{S} \widehat{p}_s \widehat{\mathbf{h}}_s \widehat{\mathbf{h}}_s^T \right) \mathbf{R}_{J_n,k}^{-1/2}$$

*in $\mathbb{R}^{J_n}$. Then, for each $j$, $\lambda_j$ is an eigenvalue of $\widetilde{L}_{n,k}^\mathcal{G}$ and $(R(\cdot, t_1), \ldots, R(\cdot, t_{J_n})) \times \mathbf{R}_{J_n,k}^{-1/2} \mathbf{u}_j$ is the corresponding eigenfunction.*

Thus, under the assumptions of Lemma 8,

(23) $$\|(R(\cdot, t_1), \ldots, R(\cdot, t_{J_n})) \mathbf{R}_{J_n,k}^{-1/2} \mathbf{u}_j - f_j\|_{\mathcal{H}_R} \xrightarrow{p} 0$$

$$\text{as } n \to \infty \text{ and then } k \to k_\infty.$$

Let

$$\widetilde{\xi}_{n,k,j} := \Psi_X^{-1}((R(\cdot,t_1),\ldots,R(\cdot,t_{J_n}))\mathbf{R}_{J_n,k}^{-1/2}\mathbf{u}_j) = (X_{t_1},\ldots,X_{t_{J_n}})\mathbf{R}_{J_n,k}^{-1/2}\mathbf{u}_j.$$

Since $\Psi_X$ is an isometric isomorphism (Proposition 2) and $\Psi(\xi_j) = f_j$, (23) is equivalent to

(24) $\qquad \|\widetilde{\xi}_{n,k,j} - \xi_j\|_{L_X^2} \xrightarrow{p} 0 \qquad$ as $n \to \infty$ and then $k \to k_\infty$.

However, since $R$ is unknown, $\widetilde{\xi}_{n,k,j}$ cannot be directly used for inference. Intuitively, $(\lambda_j, \mathbf{u}_j)$ in Lemma 9 can be estimated by the eigenvalues and eigenvectors of $\widehat{\mathbf{M}}_{n,k}$ in (16). The following result provides the justification.

LEMMA 10. *Assume that $\sup_t \mathrm{E}(X_t^4) < \infty$. Also assume that $J_n = o(n)$ and (19) holds for $m = k$, a fixed positive integer. Then*

(25) $\qquad \|\widehat{\mathbf{M}}_{n,k} - \widetilde{\mathbf{M}}_{n,k}\|_\infty \xrightarrow{p} 0 \qquad$ as $n \to \infty$.

*Also if $\mathbf{u}_j$ and $\mathbf{v}_j$ are the eigenvectors corresponding to the jth eigenvalues of $\widetilde{\mathbf{M}}_{n,k}$ and $\widehat{\mathbf{M}}_{n,k}$, respectively, we have*

(26) $\qquad \|\widehat{\xi}_{n,k,j} - \widetilde{\xi}_{n,k,j}\|_{L_X^2} \xrightarrow{p} 0 \qquad$ as $n \to \infty$.

Combining (24) and (26), we have:

THEOREM 11. *Assume that $\sup_t \mathrm{E}(X_t^4) < \infty$, $J_n = o(n)$, and either $T = \bigcup_{J=1}^\infty T_J$, or $\bigcup_{J=1}^\infty T_J$ is dense in $T$ and $R$ is continuous on $T \times T$. Also assume that (17) holds, and that (19) holds for all $m \in \mathfrak{K} = \{k_1, k_2, \ldots\}$ where $k_\ell \to k_\infty$. Then for each $j$,*

(27) $\qquad \|\widehat{\xi}_{n,k_\ell,j} - \xi_j\|_{L_X^2} \xrightarrow{p} 0 \qquad$ as $n \to \infty$ and then $\ell \to \infty$.

REMARKS. (a) The interpretation of (27) is that, under the assumptions of the theorem, there exists a sequence $\ell_n$ such that $\widehat{\xi}_{n,k_{\ell_n},j} \xrightarrow{p} \xi_j$ in $L_X^2$. In reality, $k_{\ell_n}$ is picked so that $\widehat{\mathbf{R}}_{n,J_n,k_{\ell_n}}$ and $\widehat{\mathbf{R}}_{n,J_n,k_{\ell_n}}^-$ estimate $\mathbf{R}_{J_n}$ and $\mathbf{R}_{J_n,k_{\ell_n}}^-$, respectively, well.

(b) We conjecture that the assumption $J_n = o(n)$ can be considerably relaxed. The assumption is needed because, in our proofs, we bound the distances between certain operators using the Hilbert–Schmidt norm. To relax the condition requires a different approach of bounding those distances, which is beyond our reach at this point.

**5. Numerical examples.** We now demonstrate the methodology in Section 4 with some numerical examples. Examples 1 and 2 are based on computer simulations, and Example 3 contains an analysis of real data.

In order to implement the methodology in Section 4, we need to know how to choose $k$ in $\widehat{\xi}_{n,k,j}$. Recall that $k$ is a smoothing parameter which controls the bias/variance trade-off. Also recall from the asymptotic theory [cf. (17) and (19)] that our procedure is designed to deal with situations where the effective dimension of the data is much smaller than $J_n$, the actual length of the data vector. In practice, $k$ can be chosen subjectively to ensure that $\sum_{j=1}^{k} \lambda_j(\widehat{\mathbf{R}}_{n,J_n})/\sum_{all\,j} \lambda_j(\widehat{\mathbf{R}}_{n,J_n})$ is close to 1, and yet the eigenvalues $\lambda_j(\widehat{\mathbf{R}}_{n,J_n}), 1 \leq j \leq k$, are not "too small." However, the following data-driven procedure for choosing $k$ may be useful. Consider the model

$$Y = \ell(\xi_1, \ldots, \xi_p) + \varepsilon \tag{28}$$

and assume that $\ell$ is smooth. For each feasible $k$ and $i = 1, \ldots, n$, we leave out $(\mathbf{x}_i, y_i)$ and use the rest of the data $(\mathbf{x}_{[-i]}, \mathbf{y}_{[-i]})$ to compute the $\widehat{\xi}_{n,k,j}$ in (15) and nonparametrically estimate $\ell$; use the $\widehat{\xi}_{n,k,j}$, the estimated $\ell$, and $\mathbf{x}_i$ to compute a predicted value $\widehat{y}_{i,k}$; let $CV(k) = \sum_{i=1}^{n}(y_i - \widehat{y}_{i,k})^2$, and pick $k$ to minimize $CV(k)$. Instead of leaving one datum out at a time, given enough data, we can also divide the data into training and testing samples in computing $CV$; see Example 3. These cross-validation procedures are not ideal since we need to know $p$ in advance, and the nonparametric fitting adds an extra layer of complication. A more satisfactory procedure that is free of these problems is currently not available.

The number of slices $S$ in SIR is another issue. However, it is a relatively minor one which usually does not change the outcomes of the analysis in a big way. We let $S = 10$ in all of the following examples.

EXAMPLE 1. Let $\{X(t), t \in [0,1]\}$ be a standard Brownian motion, $\varepsilon \sim N(0, 0.3^2)$, and

$$Y = \exp\left(\int_0^1 \beta(s) X(s)\,ds\right) + \varepsilon,$$

where $\beta(s) = \sin(3\pi s/2)$. Hence $\xi = \int_0^1 \beta(s) X(s)\,ds$. A sample of $n = 100$ i.i.d. $(x_i, y_i)$ were generated, where each $x_i$ was observed at 100 equally spaced time points in $[0,1]$. The first five eigenvalues of the sample covariance $\widehat{\mathbf{R}}_{n,J_n}$ are 35.17, 4.06, 1.65, 0.75 and 0.54 compared to the first five theoretical eigenvalues 0.405285, 0.045031, 0.016211, 0.008271, and 0.005003 of the Brownian motion in $L^2[0,1]$. The amounts of variation in the sample explained by the first five eigenvectors of the sample covariance cumulatively are 0.80, 0.89, 0.93, 0.94 and 0.96. The cross-validation procedure described in the beginning of this section selected $k = 2$. The plots of $\widehat{\beta}, \widehat{\xi}$ versus $\xi$, and $y$ versus $\widehat{\xi}$ are displayed in Figure 1. See (15) for the definitions of $\widehat{\beta}$ and $\widehat{\xi}$. We also estimated the link function $\ell$ by smoothing spline, which is displayed along with the plot for $y$ versus $\widehat{\xi}$. It is not surprising that $\widehat{\beta}$ is not smooth since no smoothing took place in computing it. If desired, a smoothing procedure can be incorporated in the eigendecomposition of $\widehat{\mathbf{M}}_{n,k}$ [see, e.g., Silverman (1996)]. Note that the results presented are based on one single simulation run. However, the quality of the estimates, especially for $\ell$ and $\xi$, is largely representative of what is obtained in repeated simulations. In particular, the sample correlations were seen to be averaging over 0.98 in repeated simulation runs.

EXAMPLE 2. Consider the model in which $X$ is a fractional Gaussian process on $[0,1]$ with self-similarity index $H = 0.75$ [cf. Samorodnitsky and Taqqu (1994)], and

$$Y = \tan^{-1}\left(\sum_{i=30}^{32} X_{i/121} + \sum_{i=90}^{92} X_{i/121}\right) + \varepsilon,$$

where $\varepsilon \sim N(0, 0.3^2)$. Note that, in this case, $\xi$ cannot be written as an $L^2[0,1]$ inner product of a smooth curve $\beta$ with $X$. A sample of $n = 80$ i.i.d. $(x_i, y_i)$ were generated, where each $x_i$ was observed at 120 equally spaced time points in $[0,1]$. The same methodology as in Example 1 was applied, and the results are displayed in Figure 2. The variation in the sample explained by the first four eigenvectors of the sample variance $\widehat{\mathbf{R}}_{n,J_n}$ exceeded 99%. However, cross validations picked $k = 8$. Other simulation runs produced qualitatively similar results.

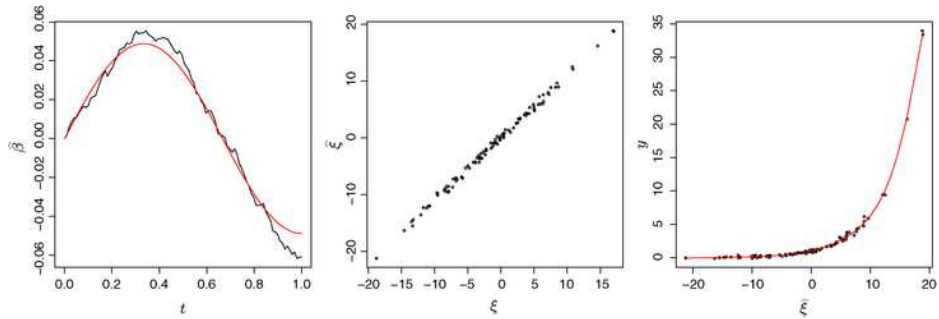

FIG. 1. *The leftmost plot is $\beta$ (smooth curve) and $\widehat{\beta}$ (nonsmooth curve) versus $t$, the middle plot is $\widehat{\xi}$ versus $\xi$, and the right plot is $y$ versus $\widehat{\xi}$.*

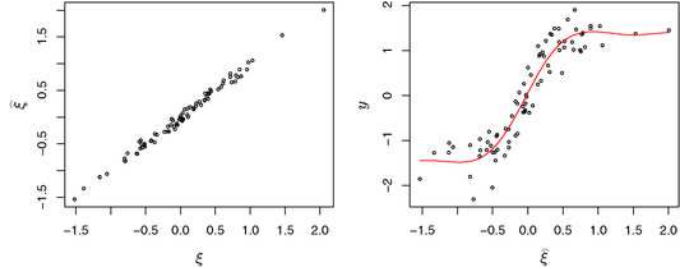

Fig. 2. *The left plot is $\widehat{\xi}$ versus $\xi$, and the right plot is $y$ versus $\widehat{\xi}$.*

Example 3. Consider a set of data recorded by the Tecator Infratec Food and Feed Analyzer, available at [http://lib.stat.cmu.edu/datasets/tecator](http://lib.stat.cmu.edu/datasets/tecator), and which were analyzed by Ferré and Yao ([2005](#)), Amato, Antoniadis and Feis ([2006](#)) and Ferraty and Vieu ([2006](#)). Each food sample contains finely minced pork meat with different contents of fat, protein and moisture. During the experiment, the spectrometer measured the spectrum of light transmitted through the sample in the region 850–1050 nanometers (nm). For each meat sample, the data consist of a 100-channel spectrum of absorption and the contents of fat, protein and moisture. The spectral data are partially observed functional data, whereas fat, protein and moisture contents are multivariate data. The spectral data are transformed to $-\log_{10}$ of their original value. In this example, we focus on the regression of spectrum $X$ on fat content $U$. In accordance with the literature, we perform the normalizing transformation $Y = \log_{10}(U/(1-U))$.

The sample size of these data is 240, and, as in Amato, Antoniadis and Feis ([2006](#)), we use the first 125 for training, and the remaining 115 for validation. The first three eigenvectors of the sample covariance $\widehat{\mathbf{R}}_{n,J_n}$ explain over 99.5% of the total variation. For different values of $k$, we used the first four of the estimated edr variables to estimate $\ell$, where the smoothing spline anova function ssanova in R [cf. Gu ([2002](#))] was our fitting algorithm. The validated prediction errors, $\{n^{-1}\sum_{i=1}^n (\widehat{y}_i - y_i)^2\}^{1/2}$, for $k = 5, 21$ and 25 were 0.06842495, 0.04481962 and 0.07414923, respectively, with $k = 21$ achieving the smallest prediction error. With $k = 21$, the two plots on the left of Figure 3 are the estimates of the first two RKHS edr functions, and the plot on the right of Figure 3 is $\widehat{y} := \widehat{\ell}(\widehat{\xi}_1, \widehat{\xi}_2, \widehat{\xi}_3, \widehat{\xi}_4)$ versus $y$ for the validation sample.

## 6. Proofs.

Proof of Proposition 2. Consider $\eta = \sum_{i=1}^m a_i X(s_i), \xi = \sum_{j=1}^n b_j X(t_j)$. By the reproducing property,

$$\langle \Psi_X(\eta), \Psi_X(\xi) \rangle_{\mathcal{H}_R} = \sum_{i=1}^m \sum_{j=1}^n a_i b_j \langle R_{s_i}, R_{t_j} \rangle_{\mathcal{H}_R} = \sum_{i=1}^m \sum_{j=1}^n a_i b_j R(s_i, t_j),$$

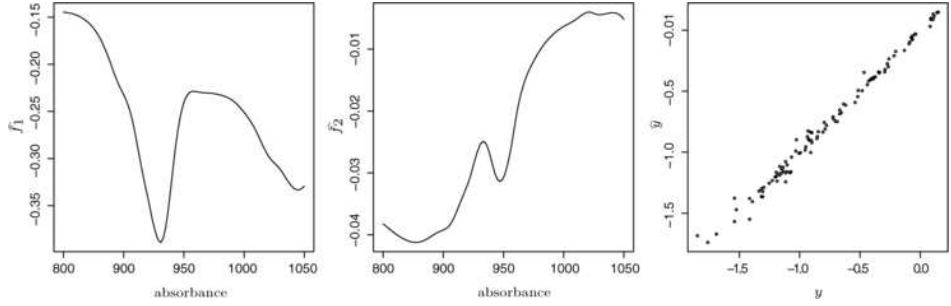

Fig. 3. *The two plots on the left describe the estimated RKHS edr functions $\widehat{f}_1$ and $\widehat{f}_2$ for Tecator data, and the plot on the right is $\widehat{y} = \widehat{\ell}(\widehat{\xi}_1, \widehat{\xi}_2, \widehat{\xi}_3, \widehat{\xi}_4)$ versus $y$.*

which is $\langle \eta, \xi \rangle_{L_X^2}$. The equality extends readily to general random variables in $L_X^2$ since random variables of the form $\eta, \xi$ are dense. $\square$

PROOF OF PROPOSITION 4. Since

$$\frac{|\sum_{i=1}^n a_i f(t_i)|^2}{\sum_{i=1}^n \sum_{j=1}^n a_i a_j K_2(t_i, t_j)} \leq \frac{|\sum_{i=1}^n a_i f(t_i)|^2}{\sum_{i=1}^n \sum_{j=1}^n a_i a_j K_1(t_i, t_j)},$$

(a) follows at once from Proposition 3. To show (b), note that for $f = \sum_{i=1}^n c_i \times K_2(\cdot, t_i)$, we have $Lf = \sum_{i=1}^n c_i K_1(\cdot, t_i) =: f_1$, and hence

$$(29) \qquad \frac{\|Lf\|_{\mathcal{H}_{K_1}}^2}{\|f\|_{\mathcal{H}_{K_2}}^2} = \frac{\|f_1\|_{\mathcal{H}_{K_1}}^2}{\|f\|_{\mathcal{H}_{K_2}}^2} = \frac{\sum_{i=1}^n \sum_{j=1}^n c_i c_j K_1(t_i, t_j)}{\sum_{i=1}^n \sum_{j=1}^n c_i c_j K_2(t_i, t_j)} \leq 1.$$

Since the set of $f$ of the above form is dense in $\mathcal{H}_{K_2}$, (29) holds for all $f \in \mathcal{H}_{K_2}$. This shows that $L$ is bounded. That $L$ is nonnegative, and self-adjoint can be seen easily by the reproducing property. $\square$

PROOF OF PROPOSITION 5. Part (a) follows at once from Fortet's formula in Proposition 3. To show (b), we focus on the case where $K$ and $f$ are continuous. Note that for any arbitrary finite set of points, $t_1, \ldots, t_n \subset T$ and constants $a_1, \ldots, a_n \in \mathbb{R}$ such that $\sum_i \sum_j a_i a_j K(t_i, t_j) \neq 0$, it follows from the continuity of $K$ and $f$, together with Fortet's formula, that

$$\frac{|\sum_{i=1}^n a_i f(t_i)|^2}{\sum_{i=1}^n \sum_{j=1}^n a_i a_j K(t_i, t_j)} \leq \lim_{n \to \infty} \|f_n\|_{\mathcal{H}_{K_n}}.$$

Then it is clear that $\lim_{n \to \infty} \|f_n\|_{\mathcal{H}_{K_n}}$ is equal to the expression on the left-hand side of (2), and (b) follows from Proposition 3. If, instead, $T$ is countable, then the above proof can be easily adapted to yield the desired conclusion and is omitted. $\square$

LEMMA 12. *Let $T$ be a separable metric space. Assume that $K_1$ and $K_2$ are covariance kernels on $T \times T$ such that:*

(a) $K_2 \gg K_1$, *and*
(b) *either $T$ is countable or $K_2$ is continuous.*

*Define a countable set $S_0 = \{s_1, s_2, \ldots\}$ which is equal to $T$ if $T$ is countable, and some arbitrary dense subset of $T$ otherwise. Denote by $L$ the dominance operator of $\mathcal{H}_{K_2}$ over $\mathcal{H}_{K_1}$, and, for $i = 1, 2$, let $K_{i,n}$ be the restriction of $K_i$ to $S_n \times S_n$ where $S_n = \{s_1, s_2, \ldots, s_n\}$. Then we can compute $\operatorname{tr}(L)$ by the formula*

$$\operatorname{tr}(L) = \lim_{n \to \infty} \operatorname{tr}(K_{1,n} K_{2,n}^-),$$

*where $K_{2,n}^-$ is the Moore–Penrose generalized inverse of $K_{2,n}$.*

PROOF. Let $K_{i,0}$ be the restriction of $K_i$ to $S_0 \times S_0$. We first establish that $\mathcal{H}_{K_{i,0}}$ and $\mathcal{H}_{K_i}$ are isometrically isomorphic. If $T = S_0$, there is nothing to prove. So we focus on the case where $K_2$ is continuous and $S_0$ is a dense subset of $T$. Note that, since $K_2 \gg K_1$, the continuity of $K_2$ implies that of $K_1$. For $s, s' \in S_0$,

$$\|K_{i,0}(\cdot, s) - K_{i,0}(\cdot, s')\|^2_{\mathcal{H}_{K_{i,0}}} = \|K_i(\cdot, s) - K_i(\cdot, s')\|^2_{\mathcal{H}_{K_i}}$$
$$= K_i(s, s) + K_i(s', s') - 2K_i(s, s'),$$

which tends to 0 if $s, s'$ both approach a fixed point $t \in T$ by continuity. By completeness $K_{i,0}(\cdot, t)|_{S_0} \in \mathcal{H}_{K_{i,0}}$ for each $t \in T$. Then it is easy to see that $\mathcal{H}_{K_{i,0}}$ and $\mathcal{H}_{K_i}$ are isometrically isomorphic. Thus, it suffices to prove

$$\operatorname{tr}(L_0) = \lim_{n \to \infty} \operatorname{tr}(K_{1,n} K_{2,n}^-),$$

where $L_0$ is the dominance operator of $\mathcal{H}_{K_{2,0}}$ over $\mathcal{H}_{K_{1,0}}$. This follows from the argument below [cf. Lukić and Beder (2001)]. Apply the Gram–Schmidt procedure to the functions $K_{2,0}(\cdot, s_i), i = 1, 2, \ldots$, to obtain a CONS $e_{1,0}, e_{2,0}, \ldots$ of $\mathcal{H}_{K_{2,0}}$. Thus,

$$\operatorname{tr}(L_0) = \sum_{i=1}^{\infty} \langle L_0 e_{i,0}, e_{i,0} \rangle_{\mathcal{H}_{K_{2,0}}}.$$

Let $L_n$ be the dominance operator of $\mathcal{H}_{K_{2,n}}$ over $\mathcal{H}_{K_{1,n}}$, and $e_{i,n}$ be the restriction of $e_{i,0}$ to $S_n$. It is clear that $e_{i,n}, 1 \le i \le n$, form an orthonormal basis for $\mathcal{H}_{K_{2,n}}$. Thus,

$$\operatorname{tr}(L_n) = \sum_{i=1}^{n} \langle L_n e_{i,n}, e_{i,n} \rangle_{\mathcal{H}_{K_{2,n}}} = \sum_{i=1}^{n} \langle L_0 e_{i,0}, e_{i,0} \rangle_{\mathcal{H}_{K_{2,0}}},$$

so that $\operatorname{tr}(L_0) = \lim_{n\to\infty} \operatorname{tr}(L_n)$. Viewing $L_n, K_{1,n}$ and $K_{2,n}$ as matrices, it follows from (b) of Proposition 4 that $L_n K_{n,2} = K_{n,1}$. Thus, $L_n = K_{1,n} K_{2,n}^-$, which completes the proof. □

LEMMA 13. *Let $T$ be a separable metric space, and let $\{U_t, t \in T\}$ be a second-order process on $T$ with mean 0 and covariance function $K_1$. Let $K_2$ be another covariance kernel on $T \times T$ such that:*

(a) $K_2 \gg K_1$, *and,*
(b) *either $T$ is countable, or both $K_2$ is continuous on $T \times T$ and the sample paths of $U$ are continuous a.s. on $T$.*

*Then* $\mathrm{P}(U \in \mathcal{H}_{K_2}) = 1$.

PROOF. Let $S_0, S_n, K_{1,n}, K_{2,n}$ and $L$ be as defined in Lemma 12. Note that $\operatorname{tr}(L) < \infty$ by the assumption $K_2 \gg K_1$. Define $U_n = U|_{S_n}$. Since $U_n$ is finite-dimensional, it is easily seen that $U_n \in \mathcal{H}_{K_{1,n}}$ a.s., which implies that $U_n \in \mathcal{H}_{K_{2,n}}$ a.s. by (a) of Proposition 4. By (7) and the property of trace,

$$\mathrm{E}(\|U_n\|_{K_{2,n}}^2) = \mathrm{E}(U_n^T K_{2,n}^- U_n) = \mathrm{E}[\operatorname{tr}(U_n^T K_{2,n}^- U_n)] = \mathrm{E}[\operatorname{tr}(U_n U_n^T K_{2,n}^-)]$$
$$= \operatorname{tr}[\mathrm{E}(U_n U_n^T) K_{2,n}^-] = \operatorname{tr}(K_{1,n} K_{2,n}^-).$$

Since $\|U_n\|_{K_{2,n}}$ is monotone by (a) of Proposition 5, it follows from the monotone convergence theorem and Lemma 12 that

$$\mathrm{E}\left[\lim_{n\to\infty} \|U_n\|_{K_{2,n}}^2\right] = \lim_{n\to\infty} \operatorname{tr}(K_{1,n} K_{2,n}^-) = \operatorname{tr}(L) < \infty. \tag{30}$$

This implies that $\lim_{n\to\infty} \|U_n\|_{K_{2,n}}^2 < \infty$ a.s., which, by (b) of Proposition 5, implies that $U \in \mathcal{H}_{K_2}$ a.s. □

PROOF OF THEOREM 6. The proof is accomplished in three steps below.
(a) Verify that $\dim(\mathcal{H}_K) \leq p$.
By definition $L_Z^2 = \overline{\operatorname{span}}\{Z_t, t \in T\}$. It follows from (IR1) and (IR2) that for each $t \in T$,

$$Z_t = \mathrm{E}(\mathrm{E}(X_t | Y, \xi_1, \ldots, \xi_p) | Y) = \mathrm{E}(\mathrm{E}(X_t | \xi_1, \ldots, \xi_p) | Y)$$
$$= \sum_{i=1}^{p} c_{i,t} \mathrm{E}(\xi_i | Y) \quad \text{a.s.}$$

for some constants $c_{i,t}$. It follows that $Z_t \in \operatorname{span}\{\mathrm{E}(\xi_i | Y), i = 1, \ldots, p\}$, $t \in T$. Consequently, $L_Z^2 \subseteq \operatorname{span}\{\mathrm{E}(\xi_i | Y), i = 1, \ldots, p\}$, and hence $\dim(\mathcal{H}_K) = \dim(L_Z^2) \leq p$.

(b) Verify that $Z \in \mathcal{H}_R$ a.s. By step (a) and (4), we conclude at once that $R \gg K$ and the dominance operator $L$ of $\mathcal{H}_R$ over $\mathcal{H}_K$ is of finite rank and hence nuclear with $\mathrm{tr}(L) < \infty$. Thus, the desired conclusion here follows from Lemma 13 under the condition (IR3).

(c) Finally, prove that $Z \in \mathcal{H}_{X,e}$ a.s. We will show that $\langle Z, h \rangle_{\mathcal{H}_R} = 0$ a.s. for any $h \in \mathcal{H}_R$ such that

$$\langle h, \Psi_X(\xi_i) \rangle_{\mathcal{H}_R} = 0, \qquad 1 \leq i \leq p. \tag{31}$$

Fix such an $h$ and let $\xi = \Psi_X^{-1}(h) \in L_X^2$. If $h = R_t$, then $\xi = X_t$, and, by the reproducing property, we obtain

$$\langle Z, h \rangle_{\mathcal{H}_R} = Z_t = \mathrm{E}(X_t|Y) = \mathrm{E}(\xi|Y).$$

Hence, in general, we have

$$\langle Z, h \rangle_{\mathcal{H}_R} = \mathrm{E}(\xi|Y) \qquad \text{for all } h \in \mathcal{H}_R.$$

By the properties of conditional expectation and (IR1),

$$\mathrm{E}(\xi|Y) = \mathrm{E}(\mathrm{E}(\xi|\xi_1, \ldots, \xi_p, Y)|Y) = \mathrm{E}(\mathrm{E}(\xi|\xi_1, \ldots, \xi_p)|Y).$$

Thus, it suffices to show that the above right-hand side equals 0, which we now do. Since by (IR2), $\mathrm{E}(\xi|\xi_1, \ldots, \xi_p) = \sum_{i=1}^p c_i \xi_i$ for some $c_i, 1 \leq i \leq p$, we have

$$\mathrm{E}(\mathrm{E}^2(\xi|\xi_1, \ldots, \xi_p)) = \mathrm{E}\left(\sum_{i=1}^p c_i \xi_i \, \mathrm{E}(\xi|\xi_1, \ldots, \xi_p)\right)$$
$$= \mathrm{E}\left(\sum_{i=1}^p c_i E(\xi \xi_i | \xi_1, \ldots, \xi_p)\right) = \sum_{i=1}^p c_i \, \mathrm{E}(\xi \xi_i),$$

which, by isometry and (31), is equal to

$$\sum_{i=1}^p c_i \langle \Psi_X(\xi), \Psi_X(\xi_i) \rangle_{\mathcal{H}_R} = \sum_{i=1}^p c_i \langle h, \Psi_X(\xi_i) \rangle_{\mathcal{H}_R} = 0.$$

Thus, $\mathrm{E}(\xi|\xi_1, \ldots, \xi_p) = 0$ a.s. and therefore $\mathrm{E}(\mathrm{E}(\xi|\xi_1, \ldots, \xi_p)|Y) = 0$ a.s. The proof is complete. $\square$

PROOF OF PROPOSITION 7. Without loss of generality, take $[a, b]$ to be $[0, 1]$ and, for convenience, let $t_0 = 0$ and $t_i = j/J, 1 \leq j \leq J$. Let $R_J$ be the discretized version of $R$:

$$R_J(s, t) = \sum_{i,j=1}^J R(t_i, t_j) I((s, t) \in [t_{i-1}, t_i) \times [t_{j-1}, t_j)).$$

Define the integral operator $Q_J : f \to \int_0^1 R_J(\cdot, y) f(y) \, dy$ on $L^2[0, 1]$. Note that $Q$ is Hilbert–Schmidt and hence has a countable number of eigenvalues. It is straightforward to verify that $Q_J$ has the same eigenvalues as $J^{-1}\mathbf{R}_J$, and $Q_J$ converges to $Q$ in uniform norm. Thus, $\lambda_j(\mathbf{R}_J) \sim J\lambda_j(Q)$ for each fixed $j$. Hence, (19) holds if the multiplicity of $\lambda_m(Q)$ is 1.

To show that (17) holds, let $\lambda_i, \phi_i$ be the eigenvalues and eigenfunctions of $Q$. By Mercer's theorem, $R(s,t) = \sum_i \lambda_i \phi_i(s) \phi_i(t)$. Define

$$R_{(k)} = \sum_{i \geq k+1} \lambda_i \phi_i(s) \phi_i(t).$$

For any $f \in \mathcal{H}_R$, write $f_{(k)} = \sum_{i \geq k+1} \langle f, \phi_i \rangle_{L^2[0,1]} \phi_i$. It is obvious that

$$\|f_{(k)}\|_{\mathcal{H}_{R_{(k)}}} \leq \|f\|_{\mathcal{H}_R}. \tag{32}$$

Now we claim that

$$\lim_{k \to \infty} \|f_{(k)}\|_{\mathcal{H}_{R_{(k)}}} = 0, \qquad f \in \mathcal{H}_R. \tag{33}$$

Given $\varepsilon > 0$, there exist some finite $M$ and constants $c_i$ such that the approximation $\tilde{f} = \sum_{i=1}^M c_m R(\cdot, t_m)$ of $f$ satisfies $\|f - \tilde{f}\|_{\mathcal{H}_R} < \varepsilon$. Write

$$\|f_{(k)}\|_{\mathcal{H}_{R_{(k)}}} = \|(f - \tilde{f} + \tilde{f})_{(k)}\|_{\mathcal{H}_{R_{(k)}}} \leq \|(f - \tilde{f})_{(k)}\|_{\mathcal{H}_{R_{(k)}}} + \|(\tilde{f})_{(k)}\|_{\mathcal{H}_{R_{(k)}}}.$$

The first term on the right-hand side is bounded by $\varepsilon$ by (32). Note that $(\tilde{f})_{(k)} = \sum_{i=1}^M c_m R_{(k)}(\cdot, t_m)$ so that

$$\|(\tilde{f})_{(k)}\|_{\mathcal{H}_{R_{(k)}}} = \mathbf{c}^T \mathbf{R}_{(k)} \mathbf{c} \to 0 \qquad \text{as } k \to \infty,$$

where $\mathbf{R}_{(k)} = \{R_{(k)}(t_i, t_j)\}_{i,j=1}^M$. This shows that

$$\limsup_{k \to \infty} \|f_{(k)}\|_{\mathcal{H}_{R_{(k)}}} < \varepsilon.$$

Since $\varepsilon$ is arbitrary, (33) follows. Now for any process $U_t, t \in T$, whose sample paths are in $\mathcal{H}_R$ a.s. and $\mathrm{E}(\|U\|_{\mathcal{H}_R}^2) < \infty$, by (32), (33) and Lebesgue's dominated convergence theorem,

$$\lim_{k \to \infty} \mathrm{E}(\|U_{(k)}\|_{\mathcal{H}_{R_{(k)}}}^2) = 0.$$

In particular,

$$\lim_{k \to \infty} \mathrm{E}(\|Z_{(k)}\|_{\mathcal{H}_{R_{(k)}}}^2) = 0,$$

where $Z = \mathrm{E}(X|Y)$. However, by the proof of Lemma 13,

$$\mathrm{E}(\|Z_{(k)}\|_{\mathcal{H}_{R_{(k)}}}^2) = \lim_{J \to \infty} \mathrm{tr}((\mathbf{R}_J^- - \mathbf{R}_{J,k}^-)\mathbf{K}_J).$$

Hence, (17) holds. □

For each $J$, let $\mathcal{H}_{R_J}$ be the subspace of $\mathcal{H}_R$ spanned by $R(\cdot, t_j), j = 1, \ldots, J$. Let $\mathbf{R}_J = \{R(t_i, t_j)\}_{i,j=1}^J$. Each $f \in \mathcal{H}_{R_J}$ can be written as $f = (R(\cdot, t_1), \ldots,$ $R(\cdot, t_J))\mathbf{c}, \mathbf{c} \in \mathbb{R}^J$, where, by the reproducing property, $\|f\|^2_{\mathcal{H}_R} = \mathbf{c}^T \mathbf{R}_J \mathbf{c}$. Thus, without loss of generality, write

$$\mathcal{H}_{R_J} = \{(R(\cdot, t_1), \ldots, R(\cdot, t_J))\mathbf{c} : \mathbf{c} \in \mathrm{Im}(\mathbf{R}_J)\}.$$

Let $\Pi_J$ be the projection operator from $\mathcal{H}_R$ into $\mathcal{H}_{R_J}$. Also define the space

$$\mathcal{H}_{R_{J,k}} = \{(R(\cdot, t_1), \ldots, R(\cdot, t_J))\mathbf{c} : \mathbf{c} \in \mathrm{Im}(\mathbf{R}_{J,k})\}$$

and the projection $\Pi_{J,k}$ from $\mathcal{H}_R$ into $\mathcal{H}_{R_{J,k}}$.

LEMMA 14. *For any $f \in \mathcal{H}_R$, and $J \geq k \geq 1$,*

(34) $$\Pi_{J,k} f = (R(\cdot, t_1), \ldots, R(\cdot, t_J))\mathbf{R}_{J,k}^- \mathbf{f},$$

*where $\mathbf{f} = (f(t_1), \ldots, f(t_J))^T$.*

PROOF. By the reproducing property, for any $\mathbf{a} \in \mathrm{Im}(\mathbf{R}_{J,k})$,

$$0 = \langle f - \Pi_{J,k} f, (R(\cdot, t_1), \ldots, R(\cdot, t_J))\mathbf{a} \rangle_{\mathcal{H}_R}$$
$$= (f(t_1), \ldots, f(t_J))\mathbf{a} - ((\Pi_{J,k} f)(t_1), \ldots, (\Pi_{J,k} f)(t_J))\mathbf{a},$$

so that

(35) $$(f(t_1), \ldots, f(t_J))\mathbf{a} = ((\Pi_{J,k} f)(t_1), \ldots, (\Pi_{J,k} f)(t_J))\mathbf{a}.$$

Write

$$\Pi_{J,k} f = (R(\cdot, t_1), \ldots, R(\cdot, t_J))\mathbf{c}, \qquad \mathbf{c} \in \mathrm{Im}(\mathbf{R}_{J,k})$$

and we will show that $\mathbf{c} = \mathbf{R}_{J,k}^- \mathbf{f}$. Evaluating both sides at $t_1, \ldots, t_J$ and pre-multiplying the resulting vectors by $\mathbf{R}_{J,k}^-$, we obtain

$$\mathbf{R}_{J,k}^- [(\Pi_{J,k} f)(t_1), \ldots, (\Pi_{J,k} f)(t_J)]^T = \mathbf{R}_{J,k}^- \mathbf{R}_J \mathbf{c} = \mathbf{c}.$$

Since the rows of $\mathbf{R}_{J,k}^-$ are in $\mathrm{Im}(\mathbf{R}_{J,k})$, it follows from (35) that

$$\mathbf{R}_{J,k}^- [(\Pi_{J,k} f)(t_1), \ldots, (\Pi_{J,k} f)(t_J)]^T = \mathbf{R}_{J,k}^- \mathbf{f}.$$

Hence, $\mathbf{c} = \mathbf{R}_{J,k}^- \mathbf{f}$ and the result follows. □

LEMMA 15. *Assume that either $T = \bigcup_{J=1}^\infty T_J$, or $\bigcup_{J=1}^\infty T_J$ is dense in $T$ and $R$ is continuous on $T \times T$. We also assume that (17) holds. Then,*

$$\mathrm{E} \|(I - \Pi_{J,k}) Z\|^2_{\mathcal{H}_R} \to 0 \qquad \text{as } J \to J_\infty \text{ and then as } k \to k_\infty,$$

*where $I$ is identity mapping.*

PROOF. If $T = \bigcup_{J=1}^{\infty} T_J$, then by definition $\mathcal{H}_R = \overline{\text{span}}\{R(\cdot, t_j), j = 1, 2, \ldots\}$. Now suppose $\bigcup_{J=1}^{\infty} T_J$ is dense in $T$, where $J_\infty = \infty$, and $R$ is continuous on $T \times T$. Then as in the proof of Lemma 12, for $t_{j_\ell} \to t$ as $\ell \to \infty$, the sequence of functions $R(\cdot, t_{j_\ell})$ is Cauchy and must converge to $R(\cdot, t)$. Hence $R(\cdot, t) \in \overline{\text{span}}\{R(\cdot, t_j), j = 1, 2, \ldots\}$ for each $t$, and we also have $\mathcal{H}_R = \overline{\text{span}}\{R(\cdot, t_j), j = 1, 2, \ldots\}$. Thus, in either case, we conclude

$$\limsup_{J \to \infty} \|(I - \Pi_J)g\|_{\mathcal{H}_R} = 0, \qquad g \in \mathcal{H}_R$$

and, by Lebesgue's dominated convergence theorem,

(36) $$\limsup_{J \to \infty} \mathrm{E}\,\|(I - \Pi_J)Z\|_{\mathcal{H}_R}^2 = 0.$$

Next we establish

(37) $$\lim_{k \to k_\infty} \limsup_{J \to \infty} \mathrm{E}\,\|(\Pi_J - \Pi_{J,k})Z\|_{\mathcal{H}_R}^2 = 0,$$

which together with (36) imply the result. By Lemma 14,

$$(\Pi_J - \Pi_{J,k})Z = (R(\cdot, t_1), \ldots, R(\cdot, t_J))(\mathbf{R}_J^- - \mathbf{R}_{J,k}^-)\mathbf{Z},$$

where $\mathbf{Z} = (Z(t_1), \ldots, Z(t_J))^T$. Hence,

$$\|(\Pi_J - \Pi_{J,k})Z\|_{\mathcal{H}_R}^2 = \mathbf{Z}^T(\mathbf{R}_J^- - \mathbf{R}_{J,k}^-)\mathbf{R}_J(\mathbf{R}_J^- - \mathbf{R}_{J,k}^-)\mathbf{Z}$$
$$= \mathrm{tr}((\mathbf{R}_J^- - \mathbf{R}_{J,k}^-)\mathbf{Z}\mathbf{Z}^T).$$

Since $\mathrm{E}(\mathbf{Z}\mathbf{Z}^T) = \mathbf{K}_J$, (37) follows from (17). □

COROLLARY 16. *Assume the conditions of Lemma 15. Let $\mathcal{F}$ be a $\sigma$-field and $Z_t^{\mathcal{F}} = \mathrm{E}(Z|\mathcal{F})$. Then*

(38) $\quad \mathrm{E}\,\|(I - \Pi_{J,k})Z^{\mathcal{F}}\|_{\mathcal{H}_R}^2 \to 0 \qquad$ *as $J \to J_\infty$ and then as $k \to k_\infty$.*

*For $h_s$ defined in (12), $s = 1, \ldots, S$,*

(39) $\quad \|(I - \Pi_{J,k})h_s\|_{\mathcal{H}_R} \to 0 \qquad$ *as $J \to J_\infty$ and then as $k \to k_\infty$.*

PROOF. Since $K^{\mathcal{F}} \leq K$, (38) follows from the same proof of the lemma with $Z^{\mathcal{F}}$ replacing $Z$ everywhere. To prove (39), letting $\mathcal{G}$ be the $\sigma$-field based on which $h_s$ is defined, we have $Z^{\mathcal{G}} = h_s$ if $Y \in I_s$. Hence,

$$\|(I - \Pi_{J,k})Z^{\mathcal{G}}\|_{\mathcal{H}_R}^2 = \sum_{s=1}^{S} \|(I - \Pi_{J,k})h_s\|_{\mathcal{H}_R}^2 I(Y \in I_s).$$

Then (39) clearly follows from this and (38). □

PROOF OF LEMMA 8. We will first show that

$$\|\widetilde{h}_s - h_s\|_{\mathcal{H}_R} \xrightarrow{p} 0 \quad \text{as } n \to \infty \text{ and then } k \to k_\infty. \tag{40}$$

Write, by Lemma 14,

$$\widetilde{h}_s - h_s = (R(\cdot, t_1), \ldots, R(\cdot, t_{J_n}))\mathbf{R}^-_{J_n,k}(\widehat{\mathbf{h}}_s - \mathbf{h}_s) + (\Pi_{J_n,k} - I)h_s.$$

The second term is taken care of by Corollary 16. To show that the first term tends to 0 in probability in $\mathcal{H}_R$, note that it is equivalent to showing that

$$(\widehat{\mathbf{h}}_s - \mathbf{h}_s)^T \mathbf{R}^-_{J_n,k}(\widehat{\mathbf{h}}_s - \mathbf{h}_s) \xrightarrow{p} 0 \quad \text{as } n \to \infty. \tag{41}$$

Let $\check{\mathbf{h}}_s = (np_s)^{-1} \sum_{i=1}^n \mathbf{X}_i I(Y_i \in I_s)$. Write

$$\mathrm{E}((\check{\mathbf{h}}_s - \mathbf{h}_s)^T \mathbf{R}^-_{J_n,k}(\check{\mathbf{h}}_s - \mathbf{h}_s)) = \mathrm{tr}[\mathbf{R}^-_{J_n,k} \mathrm{E}((\check{\mathbf{h}}_s - \mathbf{h}_s)(\check{\mathbf{h}}_s - \mathbf{h}_s)^T)].$$

By independence,

$$\mathrm{E}((\check{\mathbf{h}}_s - \mathbf{h}_s)(\check{\mathbf{h}}_s - \mathbf{h}_s)^T)$$

$$= \frac{1}{n^2 p_s^2} \sum_{i=1}^n \sum_{j=1}^n \mathrm{E}(\mathbf{X}_i \mathbf{X}_j^T I(Y_i \in I_s) I(Y_j \in I_s)) - \mathbf{h}_s \mathbf{h}_s^T$$

$$= \frac{1}{np_s^2} \mathrm{E}(\mathbf{X}_1 \mathbf{X}_1^T I(Y_1 \in I_s)) - \frac{1}{n} \mathbf{h}_s \mathbf{h}_s^T \leq \frac{1}{np_s^2} \mathbf{R}_{J_n}.$$

Thus,

$$\mathrm{E}((\check{\mathbf{h}}_s - \mathbf{h}_s)^T \mathbf{R}^-_{J_n,k}(\check{\mathbf{h}}_s - \mathbf{h}_s)) \leq \frac{1}{np_s^2} \mathrm{tr}(\mathbf{R}^-_{J_n,k} \mathbf{R}_{J_n})$$
$$= \frac{1}{np_s^2} \mathrm{tr}(\mathbf{R}^-_{J_n,k} \mathbf{R}_{J_n,k}), \tag{42}$$

which tends to 0 as $n \to \infty$, since $\mathrm{tr}(\mathbf{R}^-_{J_n,k} \mathbf{R}_{J_n,k}) \leq k$. This proves (41) with $\widehat{\mathbf{h}}_s$ replaced by $\check{\mathbf{h}}_s$. Since $\widehat{p}_s \xrightarrow{\text{a.s.}} p_s$, (41) follows as well. This completes the proof of (40).

Next for $f \in \mathcal{H}_R$ with $\|f\|_{\mathcal{H}_R} = 1$,

$$\left\| \left( \widetilde{h}_s \bigotimes_{\mathcal{H}_R} \widetilde{h}_s \right) f - \left( h_s \bigotimes_{\mathcal{H}_R} h_s \right) f \right\|_{\mathcal{H}_R}$$

$$= \|\langle \widetilde{h}_s, f \rangle_{\mathcal{H}_R} \widetilde{h}_s - \langle h_s, f \rangle_{\mathcal{H}_R} h_s\|_{\mathcal{H}_R}$$

$$\leq \|\langle \widetilde{h}_s, f \rangle_{\mathcal{H}_R} \widetilde{h}_s - \langle h_s, f \rangle_{\mathcal{H}_R} \widetilde{h}_s\|_{\mathcal{H}_R} + \|\langle h_s, f \rangle_{\mathcal{H}_R} \widetilde{h}_s - \langle h_s, f \rangle_{\mathcal{H}_R} h_s\|_{\mathcal{H}_R}$$

$$\leq \|\widetilde{h}_s - h_s\|_{\mathcal{H}_R} \|\widetilde{h}_s\|_{\mathcal{H}_R} + \|\widetilde{h}_s - h_s\|_{\mathcal{H}_R} \|h_s\|_{\mathcal{H}_R}.$$

It follows from (40) that the right-hand side tends to 0 in probability, and hence

$$\left\|\widetilde{h}_s \bigotimes_{\mathcal{H}_R} \widetilde{h}_s - h_s \bigotimes_{\mathcal{H}_R} h_s \right\|_\infty \xrightarrow{p} 0 \qquad \text{as } n \to \infty \text{ and then } k \to k_\infty.$$

The result follows from this. □

PROOF OF LEMMA 9. Consider the linear mapping $\Gamma$ that maps $\mathcal{H}_{R_{J_n,k}}$ to $\mathcal{H}_{\mathbf{R}_{J_n,k}}$ such that

$$\Gamma : (R(\cdot,t_1),\ldots,R(\cdot,t_{J_n}))\mathbf{c} \mapsto \mathbf{R}_{J_n}\mathbf{c} = \mathbf{R}_{J_n,k}\mathbf{c}, \qquad \mathbf{c} \in \text{Im}(\mathbf{R}_{J_n,k}).$$

It is easy to see that $\Gamma$ is an isometric isomorphism. The operator $\widetilde{\mathbf{L}}^{\mathcal{G}}$ in $\mathcal{H}_{\mathbf{R}_{J_n,k}}$ that corresponds to $\widetilde{L}^{\mathcal{G}}$ in $\mathcal{H}_{R_{J_n,k}}$ is

$$\widetilde{\mathbf{L}}^{\mathcal{G}} = \sum_{s=1}^{S} p_s \widetilde{\mathbf{h}}_s \bigotimes_{\mathcal{H}_{\mathbf{R}_{J_n}}} \widetilde{\mathbf{h}}_s,$$

where $\widetilde{\mathbf{h}}_s = \mathbf{P}_{J_n,k}\widehat{\mathbf{h}}_s$. Representing an eigenvector of $\widetilde{\mathbf{L}}^{\mathcal{G}}$ as $\mathbf{R}_{J_n,k}\mathbf{c}, \mathbf{c} \in \text{Im}(\mathbf{R}_{J_n,k})$, by the reproducing property, the eigenequation of $\widetilde{\mathbf{L}}^{\mathcal{G}}$ is

$$\sum_{s=1}^{S} \widehat{p}_s \mathbf{P}_{J_n,k}\widehat{\mathbf{h}}_s \widehat{\mathbf{h}}_s^T \mathbf{P}_{J_n,k}\mathbf{c} = \lambda \mathbf{R}_{J_n,k}\mathbf{c}, \qquad \mathbf{c}^T \mathbf{R}_{J_n,k}\mathbf{c} = 1, \qquad \mathbf{c} \in \text{Im}(\mathbf{R}_{J_n,k}),$$

which is equivalent to

$$\sum_{s=1}^{S} \widehat{p}_s \mathbf{R}_{J_n,k}^{-1/2} \widehat{\mathbf{h}}_s \widehat{\mathbf{h}}_s^T \mathbf{R}_{J_n,k}^{-1/2} \mathbf{u} = \lambda \mathbf{u}, \qquad \mathbf{u}^T \mathbf{u} = 1, \qquad \mathbf{c} = \mathbf{R}_{J_n,k}^{-1/2} \mathbf{u}. \qquad \square$$

For a square matrix $\mathbf{A}$ containing complex elements, let $\|A\|_{\text{HS}} = \sqrt{\text{tr}(\mathbf{A}\mathbf{A}^*)}$, where $\mathbf{A}^* = \overline{A^T}$. $\|\mathbf{A}\|_{\text{HS}}$ is known as the Hilbert–Schmidt (HS) norm or Frobenius norm of $\mathbf{A}$. See Dunford and Schwarz (1988), page 1010, or Horn and Johnson (1990), Chapter 5. Note that $\|\mathbf{A}\|_\infty \leq \|\mathbf{A}\|_{\text{HS}}$.

LEMMA 17. *Assume that $\sup_t \text{E}(X_t^4) < \infty$. Then*

$$\text{E}\|\widehat{\mathbf{R}}_{n,J_n} - \mathbf{R}_{J_n}\|_{\text{HS}}^2 \leq C J_n^2/n$$

*for some universal constant $C$.*

PROOF. By definition, $\text{E}\|\widehat{\mathbf{R}}_{n,J_n} - \mathbf{R}_{J_n}\|_{\text{HS}}^2 = \text{tr}(\text{E}(\widehat{\mathbf{R}}_{n,J_n} - \mathbf{R}_{J_n})^2)$. By independence,

$$\text{E}(\widehat{\mathbf{R}}_{n,J_n} - \mathbf{R}_{J_n})^2 = \frac{1}{n}(\text{E}[(\mathbf{X}_1\mathbf{X}_1^T)^2] - \text{E}^2(\mathbf{X}_1\mathbf{X}_1^T)) \leq \frac{1}{n}\text{E}[(\mathbf{X}_1\mathbf{X}_1^T)^2].$$

Hence,

$$\mathrm{E}\,\|\widehat{\mathbf{R}}_{n,J_n} - \mathbf{R}_{J_n}\|_{\mathrm{HS}}^2 \le \frac{1}{n}\,\mathrm{E}(\mathrm{tr}(\mathbf{X}_1\mathbf{X}_1^T))^2 = \frac{1}{n}\,\mathrm{E}(\|\mathbf{X}_1\|_{\mathbb{R}_{J_n}}^4) \le C\frac{J_n^2}{n}. \qquad \square$$

LEMMA 18. *Assume that* $\sup_t \mathrm{E}(X_t^4) < \infty$. *Also assume that* $J_n = o(n)$ *and* (19) *holds for* $m = k$, *a fixed positive integer. Then*

$$\|(\widehat{\mathbf{R}}_{n,J_n,k}^{-1/2} - \mathbf{R}_{J_n,k}^{-1/2})\mathbf{P}_{J_n,k}\|_{\mathrm{HS}} = O_p(1/\sqrt{n}).$$

PROOF. Our goal is to show that for any given $\varepsilon > 0$ there exists $\delta$ such that

$$\limsup_{n\to\infty} \mathrm{P}(\|(\widehat{\mathbf{R}}_{n,J_n,k}^{-1/2} - \mathbf{R}_{J_n,k}^{-1/2})\mathbf{P}_{J_n,k}\|_{\mathrm{HS}} > \delta/\sqrt{n}) < \varepsilon.$$

First we pick $\rho$ so that

$$\limsup_{n\to\infty} \mathrm{P}(\|\widehat{\mathbf{R}}_{n,J_n} - \mathbf{R}_{J_n}\|_{\mathrm{HS}} > \rho J_n/\sqrt{n}) < \varepsilon,$$

which is possible by Lemma 17. Below we will show that on the event

$$\|\widehat{\mathbf{R}}_{n,J_n} - \mathbf{R}_{J_n}\|_{\mathrm{HS}} \le \frac{\rho J_n}{\sqrt{n}}, \tag{43}$$

we have, for some $\delta$,

$$\|(\widehat{\mathbf{R}}_{n,J_n,k}^{-1/2} - \mathbf{R}_{J_n,k}^{-1/2})\mathbf{P}_{J_n,k}\|_{\mathrm{HS}} \le \frac{\delta}{\sqrt{n}} \qquad \text{for large } n. \tag{44}$$

Without loss of generality, assume that $\lambda_k(\mathbf{R}_J) > 0$. For definiteness, let $r$ be a constant satisfying $0 < r < \liminf_{J\to J_\infty} \frac{\rho_k(\mathbf{R}_J)}{2J}$. Denote by $i$ the imaginary unit. Let $\Lambda$ be the rectangle on the complex plane with vertices $(\lambda_1(\mathbf{R}_{J_n}) + rJ_n) - rJ_n i, (\lambda_1(\mathbf{R}_{J_n}) + rJ_n) + rJ_n i, (\lambda_k(\mathbf{R}_{J_n}) - rJ_n) + rJ_n i, (\lambda_k(\mathbf{R}_{J_n}) - rJ_n) - rJ_n i$, and let $\partial\Lambda$ be the boundary of $\Lambda$. The length $\ell(\partial\Lambda)$ of $\partial\Lambda$ is $8rJ_n + 2(\lambda_1(\mathbf{R}_J) - \lambda_k(\mathbf{R}_J))$. By (18),

$$\limsup_{n\to\infty} \frac{\ell(\partial\Lambda)}{J_n} < \infty. \tag{45}$$

Since $\rho/\sqrt{n} \to 0$, by Corollary 4 on page 1090 of Dunford and Schwarz (1988) and (43), we have, for large $n$ and uniformly for all $j$,

$$|\lambda_j(\widehat{\mathbf{R}}_{n,J_n}) - \lambda_j(\mathbf{R}_{J_n})| \le \|\widehat{\mathbf{R}}_{n,J_n} - \mathbf{R}_{J_n}\|_\infty \le \|\widehat{\mathbf{R}}_{n,J_n} - \mathbf{R}_{J_n}\|_{\mathrm{HS}} < rJ_n.$$

Thus, $\Lambda$ contains $\lambda_j(\mathbf{R}_{J_n})$ and $\lambda_j(\widehat{\mathbf{R}}_{n,J_n})$, but no other eigenvalues of either $\mathbf{R}_{J_n}$ or $\widehat{\mathbf{R}}_{n,J_n}$. Also $\Lambda$ does not contain the complex origin. Let

$$\mathfrak{R}_{J_n}(z) = (z - \mathbf{R}_{J_n})^{-1} \quad \text{and} \quad \widehat{\mathfrak{R}}_{n,J_n}(z) = (z - \widehat{\mathbf{R}}_{n,J_n})^{-1}$$

be the resolvent of $\mathbf{R}_{J_n}$ and $\widehat{\mathbf{R}}_{n,J_n}$, respectively, where $z$ is complex argument restricted to the respective resolvent sets. By the Cauchy integral formula [cf. Dunford and Schwarz (1988), page 568],

$$\widehat{\mathbf{R}}_{n,J_n,k}^{-1/2} - \mathbf{R}_{J_n,k}^{-1/2} = \frac{1}{2\pi i} \oint_{\partial \Lambda} z^{-1/2} [\widehat{\mathfrak{R}}_{n,J_n}(z) - \mathfrak{R}_{J_n}(z)] \, dz$$

and hence

(46)
$$\|(\widehat{\mathbf{R}}_{n,J_n,k}^{-1/2} - \mathbf{R}_{J_n,k}^{-1/2})\mathbf{P}_{J_n,k}\|_{\mathrm{HS}}$$
$$\leq \frac{1}{2\pi} \oint_{\partial \Lambda} |z|^{-1/2} \cdot \|(\widehat{\mathfrak{R}}_{n,J_n}(z) - \mathfrak{R}_{J_n}(z))\mathbf{P}_{J_n,k}\|_{\mathrm{HS}} \, dz.$$

Note that $\mathfrak{R}_{J_n}(z)\mathbf{P}_{J_n,k} = (z - \mathbf{R}_{J_n,k})^{-1}$, which is the resolvent of $\mathbf{R}_{J_n,k}$ and we denote it subsequently as $\mathfrak{R}_{J_n,k}(z)$. Observe that

(47) $$\sup_{z\in\Lambda} \|\mathfrak{R}_{J_n}(z)\|_{\mathrm{HS}}^2 = \sup_{z\in\Lambda} \sum_{i=1}^{J_n} |z - \lambda_i(\mathbf{R}_{J_n})|^{-2} \leq \frac{J_n}{(rJ_n)^2} = \frac{1}{r^2 J_n}$$

and

(48) $$\sup_{z\in\Lambda} \|\mathfrak{R}_{J_n,k}(z)\|_{\mathrm{HS}}^2 = \sup_{z\in\Lambda} \sum_{i=1}^{k} |z - \lambda_i(\mathbf{R}_{J_n})|^{-2} \leq \frac{k}{(rJ_n)^2}.$$

Write

(49)
$$\widehat{\mathfrak{R}}_{n,J_n}(z)\mathbf{P}_{J_n,k} = (z - \mathbf{R}_{J_n} - \widehat{\mathbf{R}}_{n,J_n} + \mathbf{R}_{J_n})^{-1}\mathbf{P}_{J_n,k}$$
$$= (\mathfrak{R}_{J_n}(z)^{-1} - \widehat{\mathbf{R}}_{n,J_n} + \mathbf{R}_{J_n})^{-1}\mathbf{P}_{J_n,k}$$
$$= (\mathbf{I} - \mathfrak{R}_{J_n}(z)(\widehat{\mathbf{R}}_{n,J_n} - \mathbf{R}_{J_n}))^{-1}\mathfrak{R}_{J_n,k}(z).$$

By (43), (47), the fact that $\|\mathbf{AB}\|_{\mathrm{HS}} \leq \|\mathbf{A}\|_{\mathrm{HS}}\|\mathbf{B}\|_{\mathrm{HS}}$, and the assumption $J_n = o(n)$,

(50) $$\sup_{z\in\Lambda} \|\mathfrak{R}_{J_n}(z)(\widehat{\mathbf{R}}_{n,J_n} - \mathbf{R}_{J_n})\|_{\mathrm{HS}} < \frac{\rho}{r}\sqrt{\frac{J_n}{n}} < 1 \qquad \text{for large } n.$$

A standard argument [cf. (3.3) of Gohberg and Kreĭn (1969)] shows that

(51) $$(\mathbf{I} - \mathfrak{R}_{J_n}(z)(\widehat{\mathbf{R}}_{n,J_n} - \mathbf{R}_{J_n}))^{-1} = \sum_{j\geq 0} [\mathfrak{R}_{J_n}(z)(\widehat{\mathbf{R}}_{n,J_n} - \mathbf{R}_{J_n})]^j.$$

By (49) and (51),

$$(\widehat{\mathfrak{R}}_{n,J_n}(z) - \mathfrak{R}_{J_n}(z))\mathbf{P}_{J_n,k} = \sum_{j\geq 1} [\mathfrak{R}_{J_n}(z)(\widehat{\mathbf{R}}_{n,J_n} - \mathbf{R}_{J_n})]^j \mathfrak{R}_{J_n,k}(z)$$

and by the triangle inequality and (50),

$$
\begin{aligned}
&\|(\widehat{\mathfrak{R}}_{n,J_n}(z) - \mathfrak{R}_{J_n}(z))\mathbf{P}_{J_n,k}\|_{\mathrm{HS}} \\
&\qquad \leq \frac{\|\mathfrak{R}_{J_n}(z)\|_{\mathrm{HS}}\|\widehat{\mathbf{R}}_{n,J_n} - \mathbf{R}_{J_n}\|_{\mathrm{HS}}\|\mathfrak{R}_{J_n,k}(z)\|_{\mathrm{HS}}}{1 - \|\mathfrak{R}_{J_n}(z)(\widehat{\mathbf{R}}_{n,J_n} - \mathbf{R}_{J_n})\|_{\mathrm{HS}}}.
\end{aligned}
\tag{52}
$$

Finally,

$$
\begin{aligned}
\sup_{z\in\partial\Lambda} |z|^{-1/2} &= (\lambda_k(\mathbf{R}_{J_n}) - rJ_n)^{-1/2} \leq (\lambda_{k+1}(\mathbf{R}_{J_n}) + rJ_n)^{-1/2} \\
&< (rJ_n)^{-1/2}.
\end{aligned}
\tag{53}
$$

By (46) and (52),

$$
\begin{aligned}
&\|(\widehat{\mathbf{R}}_{n,J_n,k}^{-1/2} - \mathbf{R}_{J_n,k}^{-1/2})\mathbf{P}_{J_n,k}\|_{\mathrm{HS}} \\
&\qquad \leq \frac{\ell(\partial\Lambda)}{2\pi} \sup_{z\in\partial\Lambda} \frac{|z|^{-1/2}\cdot\|\mathfrak{R}_{J_n}(z)\|_{\mathrm{HS}}\|\widehat{\mathbf{R}}_{n,J_n} - \mathbf{R}_{J_n}\|_{\mathrm{HS}}\|\mathfrak{R}_{J_n,k}(z)\|_{\mathrm{HS}}}{1 - \|\mathfrak{R}_{J_n}(z)(\widehat{\mathbf{R}}_{n,J_n} - \mathbf{R}_{J_n})\|_{\mathrm{HS}}},
\end{aligned}
$$

from which (44) using (43), (45), (47), (48), (50) and (53). □

PROOF OF LEMMA 10. Write

$$
\begin{aligned}
\|\widehat{\mathbf{R}}_{n,J_n,k}^{-1/2}\widehat{\mathbf{h}}_s - \mathbf{R}_{J_n,k}^{-1/2}\widehat{\mathbf{h}}_s\|_{\mathbb{R}^{J_n}}^2 &= \widehat{\mathbf{h}}_s^T(\widehat{\mathbf{R}}_{n,J_n,k}^{-1/2} - \mathbf{R}_{J_n,k}^{-1/2})^2\widehat{\mathbf{h}}_s \\
&= \mathrm{tr}((\widehat{\mathbf{R}}_{n,J_n,k}^{-1/2} - \mathbf{R}_{J_n,k}^{-1/2})^2(\widehat{\mathbf{h}}_s\widehat{\mathbf{h}}_s^T)) \\
&= \|(\widehat{\mathbf{R}}_{n,J_n,k}^{-1/2} - \mathbf{R}_{J_n,k}^{-1/2})(\widehat{\mathbf{h}}_s\widehat{\mathbf{h}}_s^T)^{1/2}\|_{\mathrm{HS}}^2 \\
&= \|(\widehat{\mathbf{R}}_{n,J_n,k}^{-1/2}\mathbf{R}_{J_n,k}^{1/2} - \mathbf{P}_{J_n,k})\mathbf{R}_{J_n,k}^{-1/2}(\widehat{\mathbf{h}}_s\widehat{\mathbf{h}}_s^T)^{1/2}\|_{\mathrm{HS}}^2.
\end{aligned}
$$

Since $\|\mathbf{AB}\|_{\mathrm{HS}} \leq \|\mathbf{A}\|_{\mathrm{HS}}\|\mathbf{B}\|_{\mathrm{HS}}$, we conclude

$$
\begin{aligned}
&\|\widehat{\mathbf{R}}_{n,J_n,k}^{-1/2}\widehat{\mathbf{h}}_s - \mathbf{R}_{J_n,k}^{-1/2}\widehat{\mathbf{h}}_s\|_{\mathbb{R}^{J_n}} \\
&\qquad \leq \|(\widehat{\mathbf{R}}_{n,J_n,k}^{-1/2}\mathbf{R}_{J_n,k}^{1/2} - \mathbf{P}_{J_n,k})\|_{\mathrm{HS}}\|\mathbf{R}_{J_n,k}^{-1/2}(\widehat{\mathbf{h}}_s\widehat{\mathbf{h}}_s^T)^{1/2}\|_{\mathrm{HS}}.
\end{aligned}
\tag{54}
$$

We first address the second term of the right-hand side of (54). By definition and the same argument that leads to (42),

$$
\|\mathbf{R}_{J_n,k}^{-1/2}(\widehat{\mathbf{h}}_s\widehat{\mathbf{h}}_s^T)^{1/2}\|_{\mathrm{HS}}^2 = \mathrm{tr}(\mathbf{R}_{J_n,k}^{-}\widehat{\mathbf{h}}_s\widehat{\mathbf{h}}_s^T) = O_p(1).
\tag{55}
$$

We next deal with the first terms on the right-hand side of (54). Write

$$
\widehat{\mathbf{R}}_{n,J_n,k}^{-1/2}\mathbf{R}_{J_n,k}^{1/2} - \mathbf{P}_{J_n,k} = (\widehat{\mathbf{R}}_{n,J_n,k}^{-1/2} - \mathbf{R}_{J_n,k}^{-1/2})\mathbf{R}_{J_n,k}^{1/2}.
$$

It follows from Lemma 18 that

$$\|\widehat{\mathbf{R}}_{n,J_n,k}^{-1/2}\mathbf{R}_{J_n,k}^{1/2} - \mathbf{P}_{J_n,k}\|_{\text{HS}} \le \|(\widehat{\mathbf{R}}_{n,J_n,k}^{-1/2} - \mathbf{R}_{J_n,k}^{-1/2})\mathbf{P}_{J_n,k}\|_{\text{HS}}\|\mathbf{R}_{J_n,k}^{1/2}\|_{\text{HS}}$$
$$= O_p(\sqrt{J_n/n}) = o_p(1). \quad (56)$$

By (54)–(56), we conclude that (25) holds.

Next,

$$\|\widehat{\widetilde{\xi}}_{n,k,j} - \widetilde{\xi}_{n,k,j}\|_{L_X^2}^2$$
$$= (\mathbf{R}_{J_n,k}^{-1/2}\mathbf{u}_j - \widehat{\mathbf{R}}_{n,J_n,k}^{-1/2}\mathbf{v}_j)^T \mathbf{R}_{J_n,k}(\mathbf{R}_{J_n,k}^{-1/2}\mathbf{u}_j - \widehat{\mathbf{R}}_{n,J_n,k}^{-1/2}\mathbf{v}_j)$$
$$= (\mathbf{R}_{J_n,k}^{-1/2}(\mathbf{u}_j - \mathbf{v}_j) - (\widehat{\mathbf{R}}_{n,J_n,k}^{-1/2} - \mathbf{R}_{J_n,k}^{-1/2})\mathbf{v}_j)^T \mathbf{R}_{J_n,k}$$
$$\quad \times (\mathbf{R}_{J_n,k}^{-1/2}(\mathbf{u}_j - \mathbf{v}_j) - (\widehat{\mathbf{R}}_{n,J_n,k}^{-1/2} - \mathbf{R}_{J_n,k}^{-1/2})\mathbf{v}_j)$$
$$\le 2(\mathbf{u}_j - \mathbf{v}_j)^T \mathbf{P}_{J_n,k}(\mathbf{u}_j - \mathbf{v}_j)$$
$$\quad + 2\mathbf{v}^T(\widehat{\mathbf{R}}_{n,J_n,k}^{-1/2} - \mathbf{R}_{J_n,k}^{-1/2})\mathbf{R}_{J_n,k}(\widehat{\mathbf{R}}_{n,J_n,k}^{-1/2} - \mathbf{R}_{J_n,k}^{-1/2})\mathbf{v}$$
$$\le 2\|\mathbf{u}_j - \mathbf{v}_j\|_{\mathbb{R}^{J_n}}^2 + 2\|(\widehat{\mathbf{R}}_{n,J_n,k}^{-1/2} - \mathbf{R}_{J_n,k}^{-1/2})\mathbf{P}_{J,k}\|_{\text{HS}}^2\|\mathbf{R}_{J_n,k}^{1/2}\|_{\text{HS}}^2.$$

The first term tends to 0 in probability by the previous part, (25), whereas the second term converges to 0 in probability as in (56). $\square$

**Acknowledgment.** We are deeply indebted to the referees and Associate Editor for their many helpful comments, without which the present version of the paper would not be possible.

Department of Statistics  
University of Michigan  
439 West Hall  
1085 South University  
Ann Arbor, Michigan 48109  
USA  
E-mail: thsing@umich.edu

Clinical Biostatistics  
Sanofi-Aventis  
200 Crossings Blvd.  
P.O. Box 6890  
Bridgewater, New Jersey 08807  
USA  
E-mail: haobo.ren@sanofi-aventis.com